\documentclass{article}

\usepackage{bm}
\usepackage{color}
\usepackage{stmaryrd}
\usepackage[hang]{subfigure}
\RequirePackage[%
  papersize={15cm,24cm},
 hmargin=1cm,%
  vmargin=1.6cm,%
  head=0.2cm,
  headsep=0pt,%
  foot=0.8cm
  ]{geometry}

\newcommand\bbR{\mathbb{R}}
\newcommand\bbN{\mathbb{N}}

\newcommand\dd{\mathrm{d}}

\newcommand\bx{\bm{x}}
\newcommand\bz{\bm{z}}
\newcommand\bF{\bm{F}}
\newcommand\bG{\bm{G}}
\newcommand\bu{\bm{u}}

\newcommand\bw{\bm{w}}
\newcommand\bV{\bm{V}}
\newcommand\bU{\bm{U}}

\newcommand\xl{{i-\frac12}}
\newcommand\xr{{i+\frac12}}
\newcommand\yl{{j-\frac12}}
\newcommand\yr{{j+\frac12}}

\usepackage{amsmath}
\usepackage{amsfonts}
 \usepackage{graphicx}
\newcommand\pd[2]{\dfrac{\partial {#1}}{\partial {#2}}}

\newcommand\abs[1]{\lvert #1 \rvert}
\newcommand\jump[1]{\llbracket #1 \rrbracket}
\newcommand\mean[1]{\{\!\!\{ #1 \}\!\!\}}
\newcommand\meanln[1]{\{\!\!\{ #1 \}\!\!\}^{\text{ln}}}
\newcommand\meanL[1]{\{\!\!\{ #1 \}\!\!\}^{\text{Lor}}}
\newcommand\meanLx[1]{\{\!\!\{ #1 \}\!\!\}^{\text{Lorx}}}
\newcommand\meanLy[1]{\{\!\!\{ #1 \}\!\!\}^{\text{Lory}}}
\newcommand\jumpangle[1]{\langle\!\langle #1 \rangle\!\rangle}

\newtheorem{theorem}{Theorem}
\newtheorem{remark}{Remark}
\newtheorem{example}{Example}
\newtheorem{definition}{Definition}

\begin{document}


\title{High-order accurate entropy stable finite difference schemes for one- and two-dimensional special relativistic hydrodynamics}

\author{Junming Duan and Huazhong Tang\\
 HEDPS, CAPT \& LMAM, School of Mathematical Sciences, \\
 Peking University,
Beijing 100871,  China}

\maketitle
\begin{abstract}
This paper develops the high-order accurate entropy stable  finite difference schemes for
one- and two-dimensional special relativistic hydrodynamic equations. The schemes
are built on the entropy conservative flux and the weighted essentially non-oscillatory (WENO)
technique as well as explicit Runge-Kutta time discretization.
The key is  to technically construct the { affordable} entropy conservative flux of the  semi-discrete second-order accurate entropy conservative schemes satisfying the semi-discrete entropy equality for the found convex entropy pair.
As soon as  the entropy conservative flux  is derived,
the dissipation term can be added to give the semi-discrete entropy stable schemes satisfying the
semi-discrete entropy inequality with the given convex entropy function.
The WENO reconstruction  for the scaled entropy variables
and the high-order explicit Runge-Kutta time discretization
{
are implemented to obtain the fully-discrete high-order entropy stable schemes.
}
Several numerical tests are conducted to validate the
accuracy and the ability to capture discontinuities  of our
entropy stable schemes.



\end{abstract}


\section{Introduction}
This paper is concerned with the high-order accurate numerical schemes for the one- and two-dimensional special relativistic hydrodynamic (RHD) equations, which in the laboratory frame,
can be cast in the  divergence form 
\begin{equation}\label{eq:RHD}
  \pd{\bU}{t}+\sum_{\ell=1}^{d}\pd{\bF_\ell(\bU)}{x_\ell}=0,
\end{equation}
where $\bU$ and $\bF_\ell$ are respectively
the conservative vector and the flux vector in the $x_\ell$-direction
and defined by
\begin{align}\label{eq:RHD2}\begin{aligned}
  &\bU=(D,m_1,\cdots,m_d,E)^\mathrm{T},\\
  &\bF_\ell=(Du_\ell,m_1u_\ell+p\delta_{1,\ell},\cdots,
  m_du_\ell+p\delta_{d,\ell},m_\ell)^\mathrm{T},\quad \ell=1,\ldots,d,
\end{aligned}\end{align}
with the mass density $D=\rho W$, the momentum density
$\bm{m}=(m_1,\cdots,m_d)^\mathrm{T}=DhW\bu$, and
the energy density $E=DhW-p$.
Here, $d=1$ or 2, $\rho,p$ and $\bu=(u_1,\cdots,u_d)^\mathrm{T}$ denote the rest-mass
density, the kinetic pressure, and the fluid velocity, respectively. Moreover,
$W=1/\sqrt{1-(u_1^2+\cdots+u_d^2)}$ is the Lorentz factor and $h$
is the specific enthalpy defined by $h=1+e+p/\rho$ with units in which the speed
of light is equal to one, and  the specific internal energy $e$.
The system \eqref{eq:RHD}-\eqref{eq:RHD2} should be closed by using the equation of state (EOS). This paper
 will only consider the ideal-fluid EOS
\begin{equation*}
  p=(\Gamma-1)\rho e,
\end{equation*}
with the adiabatic index $\Gamma\in(1,2]$.
Because there is no explicit expression for the primitive variables
$(\rho,\bu^\mathrm{T},p)$ and the flux $\bF_\ell$ in terms of $\bU$,
in order to recover the values of the primitive variables and the flux from the given $\bU$,
a nonlinear algebraic equation such as
\begin{equation*}
  E+p=DW+\dfrac{\Gamma}{\Gamma-1}pW^2,
\end{equation*}
has to be numerically solved to obtain the pressure $p$,  and
then the rest-mass density $\rho$, the specific enthalpy
$h$, and the velocity $\bu$ can be orderly calculated by
\begin{equation*}
  \rho=\dfrac{D}{W},\quad h=1+\dfrac{\Gamma p}{(\Gamma-1)\rho},\quad
  \bu=\dfrac{\bm{m}}{Dh}.
\end{equation*}
The relativistic description for the dynamics of the fluid
(gas) at nearly speed of light should be considered when  we investigate the astrophysical
 phenomena  from stellar to galactic scales, e.g. coalescing neutron stars, core
collapse supernovae,  active galactic nuclei, superluminal jets, formation of black holes,
and gamma-ray bursts etc.
Due to the relativistic effect, 
the nonlinearity of the system \eqref{eq:RHD}-\eqref{eq:RHD2} becomes much stronger than the non-relativistic
case so that its analytic treatment is extremely difficult and challenging.
Numerical simulation is a primary   way to help our understanding
the physical mechanisms in the RHD.
It can be traced back to  the artificial viscosity method
for the  RHD equations in the Lagrangian coordinates \cite{may1,may2}
and in the Eulerian coordinates \cite{wilson}. 
After those, the modern shock-captured methods for the RHD equations
would not be noticed until later. Here listed are some works:
the Harten-Lax-van Leer method \cite{schneider},
the two-shock solvers \cite{balsara,dai},
the Roe solver \cite{eulderink},
the essentially non-oscillatory (ENO) and the weighted ENO (WENO) methods \cite{dolezal,zanna},
the piecewise parabolic methods \cite{marti3,mignone2},
the adaptive mesh refinement method \cite{zw},
the Runge-Kutta discontinuous Galerkin methods with WENO limiter \cite{zhao},
the direct Eulerian generalized Riemann problem schemes
\cite{yz1,yz2,wu2014,wu2016}, the adaptive moving mesh methods \cite{he1,he2}, and so on.
The readers are  referred to the early review articles \cite{marti2,font,marti2015}
and references therein.
Recently, the properties of the admissible state set and the physical-constraints-preserving (PCP) numerical schemes were  well studied for the RHD, see
\cite{wu2015,wu2017,wu2017a} and \cite{qin},
and for the special relativistic
magnetohydrodynamics \cite{wu2017b,wu2018}. The  PCP schemes
satisfy that both the rest-mass density and the kinetic pressure  are  positive and
the magnitude of the fluid velocity is less than the speed of light.
Motivated by \cite{wu2017b,wu2018}, the positivity-preserving  schemes
 for the non-relativistic ideal magnetohydrodynamics were successfully studied in
\cite{wu2018MHDa,wu2018MHDb}.

It is well known that the weak solution of the quasi-linear hyperbolic conservation laws
might not be unique so that the  entropy conditions are needed
to single out the physical relevant solution among all weak solutions.
\begin{definition}[Entropy function]
  A strictly convex function $\eta(\bU)$ is called an entropy function for the
  system \eqref{eq:RHD}-\eqref{eq:RHD2} if there are associated entropy fluxes $q_\ell(\bU)$
  such that
  \begin{equation}\label{eq:entropy}
    q_\ell'(\bU)=\bV^\mathrm{T}\bF_\ell'(\bU),\quad \ell=1,\cdots,d,
  \end{equation}
  where $\bV=\eta'(\bU)^\mathrm{T}$ is called the entropy variables, and
  $(\eta,q_\ell)$ is an entropy pair.
\end{definition}
For the smooth solutions of \eqref{eq:RHD}-\eqref{eq:RHD2},
multiplying \eqref{eq:RHD} by $\bV^\mathrm{T}$ gives the entropy identity
\begin{equation*}
  \eta(\bU)_t+\sum_{l=1}^{d} q_l(\bU)_{x_l} = 0.
\end{equation*}
However, if the solutions contain discontinuity, then the above identity
does not hold.
\begin{definition}[Entropy solution]
  A weak solution $\bU$ of \eqref{eq:RHD} is called an entropy solution if for
  all entropy functions $\eta$, the inequality
  \begin{equation}\label{eq:entropyineq}
    \eta(\bU)_t+\sum_{l=1}^{d} q_l(\bU)_{x_l} \leq 0,
  \end{equation}
  holds in the sense of distributions.
\end{definition}
Formally, integrating \eqref{eq:entropyineq} in space and imposing a periodic or
no-inflow boundary conditions gets the inequality $\dfrac{\dd}{\dd t}\int_{\bbR}\eta(\bU)\dd \bx\leq 0$,
which can be converted into an a priori
estimate on the solution of \eqref{eq:RHD} in a suitable $L^p$ space
if $\eta$ is convex \cite{Fjordholm2012,Dafermos}.
The entropy conditions  are of great importance
 in the well-posedness of hyperbolic conservation laws,
 thus it is reasonable  to seek the entropy stable schemes of \eqref{eq:RHD}, satisfying a discrete or semi-discrete version of the entropy inequality \eqref{eq:entropyineq}.
For the smooth solutions of the special RHD equations \eqref{eq:RHD}-\eqref{eq:RHD2},  the thermodynamic entropy
\begin{equation*}
  S=\ln(p)-\Gamma\ln(\rho),
\end{equation*}
satisfies
\begin{equation*}
  \pd{(\rho WS)}{t}+\sum_{\ell=1}^{d} \pd{(\rho u_\ell WS)}{x_\ell}=0,
\end{equation*}
thus  an entropy pair of \eqref{eq:RHD}-\eqref{eq:RHD2} can be defined by
\begin{equation*}
  \eta(\bU)=\dfrac{-\rho WS}{\Gamma-1},\quad q_\ell(\bU)=\dfrac{-\rho u_\ell WS}{\Gamma-1},\quad
  \ell=1,\cdots,d,
\end{equation*}
and corresponding entropy variables $\bV=\eta'(\bU)^\mathrm{T}$ can
be explicitly given by
\begin{equation*}
  \bV=\left(\dfrac{\Gamma-S}{\Gamma-1}+\dfrac{\rho}{p},
  \dfrac{\rho  W\bu^\mathrm{T}}{p},  -\dfrac{\rho W}{p} \right)^\mathrm{T},
\end{equation*}
which gives the ``potential''  $\psi_\ell:=\bV^T \bF_\ell(\bU) -q_\ell(\bU)=\rho Wu_\ell$, $\ell=1,\cdots,d$.
It can be verified that, for $d=1,2$,
the matrices $\pd{\bU}{\bV}$ and $\pd{\bF_\ell}{\bU}\pd{\bU}{\bV}$ are symmetric and $\pd{\bU}{\bV}$ is
positive definite so that \eqref{eq:RHD}-\eqref{eq:RHD2} can be symmetrized
with the above entropy pair. 
{Such symmetrization of the RHD equations \eqref{eq:RHD}-\eqref{eq:RHD2} will be useful in deriving
a set of particular scaled eigenvectors used for  designing the dissipation term of
  the entropy stable schemes in Section \ref{sec:2.2}.
}

For the scalar conservation laws, the conservative monotone schemes
were shown that they were nonlinearly stable and satisfied the discrete entropy conditions
so that they could converge to the entropy solution \cite{Harten,Crandall1980}.
A class of so-called E-schemes satisfying the entropy condition for any convex entropy were studied
in \cite{Osher1984,Osher1988}. Those schemes are only first-order accurate.
Due to the fact that it is  basically impossible to show that the high-order schemes of
the scalar conservation laws  and the schemes for the system of hyperbolic conservation laws
satisfy the entropy inequality for any convex entropy pair,
the researcher tries to study the high-order accurate
entropy stable schemes, which satisfy the entropy inequality for a given entropy pair.
The second-order entropy conservative schemes (satisfying the discrete entropy identity) were built in
\cite{Tadmor1987,Tadmor2004}, and their higher-order extension was introduced in \cite{Lefloch2002}.
It is known that an entropy conservative scheme may become oscillatory near the
shock wave, thus the additional dissipation terms has to be added to the entropy
conservative schemes to obtain the entropy stable schemes.
Combining the entropy conservative flux of the entropy conservative schemes
with the ``sign'' property of the ENO reconstructions, the
arbitrary high-order entropy stable schemes were constructed by using high-order dissipation terms \cite{Fjordholm2012}.
Some entropy stable schemes based on the DG framework were
studied, such as \cite{Barth1999,Hiltebrand2014} in the space-time DG formulation, the entropy
stable nodal DG schemes using suitable quadrature rules \cite{Chen2017}, and its
extension to magnetohydrodynamics equations \cite{Liu2018}.
The  entropy stable schemes based on
summation-by-parts (SBP) operators were developed for the Navier-Stokes equations
\cite{Fisher2013}.
The existing works
do not address the entropy stable
schemes for the special or general RHD equations.

The paper aims at constructing the high-order accurate entropy stable schemes for one- and
two-dimensional special RHD equations \eqref{eq:RHD}-\eqref{eq:RHD2}.
The key is to  technically construct the {affordable} entropy conservative flux of the
semi-discrete second-order accurate entropy conservative schemes satisfying the
semi-discrete entropy equality for the found convex entropy pair.
The paper is organized as follows.
Section \ref{section:1D}
introduces the entropy conservative fluxes, entropy conservative schemes,
and entropy stable schemes for the  one-dimensional special RHD equations.
Section \ref{section:2D} introduces our schemes for the two-dimensional special RHD equations.
Several one- and two-dimensional numerical tests are
conducted in Section \ref{section:Num} to validate the effectiveness of our schemes.
Some conclusions are summarized in Section \ref{section:conclusion}.

%

\section{One-dimensional schemes}\label{section:1D}
This section considers the one-dimensional special RHD equations, i.e.,
\eqref{eq:RHD}-\eqref{eq:RHD2} with $d=1$. For the sake of convenience,
the notations $\bF_1,u_1,m_1$ and $x_1$ are
replaced with $\bF,u,m$ and $x$, respectively,
so that the flux, the  entropy pair and the potential become
 $\bF=(Du,mu+p,m)^\mathrm{T}$,
$\eta=\dfrac{-\rho WS}{\Gamma-1}$,  $q=\dfrac{-\rho u WS}{\Gamma-1}$, and
$\psi=\rho uW$,  respectively.

Let us consider a uniform mesh $x_1<x_2<\cdots<x_N$ with the step size $\Delta
x=x_i-x_{i-1},~i=2,\cdots,N$ and the semi-discrete conservative finite
difference scheme
\begin{equation}\label{eq:semi}
  \dfrac{\dd}{\dd t}\bU_i(t)=-\dfrac{1}{\Delta
  x}\left(\hat{\bF}_{\xr}(t)-\hat{\bF}_{\xl}(t)\right),
\end{equation}
where $\bU_i(t)$ approximates the point value of $\bU(x_i,t)$ and
$\hat{\bF}_{\xr}$ is the numerical flux approximating $\bF$ at $x_{\xr}=x_i+ {\Delta x}/{2}$.

\subsection{Entropy conservative fluxes}
\begin{definition}[Entropy conservative scheme]
  The scheme \eqref{eq:semi} is called entropy conservative scheme if its solution
  satisfies a semi-discrete entropy equality
  \begin{equation}\label{eq:discEnEq}
    \dfrac{\dd}{\dd t}\eta(\bU_i(t))=-\dfrac{1}{\Delta
    x}\left(\tilde{q}_{\xr}(t)-\tilde{q}_{\xl}(t)\right),
  \end{equation}
  for some numerical entropy flux $\tilde{q}_{\xr}$ consistent with $q$.
\end{definition}
\begin{theorem}[Tadmor\cite{Tadmor1987}]
  If  a consistent numerical flux $\tilde{\bm{F}}_{\xr}$ satisfies
  \begin{equation}\label{eq:conserFlux}
    \jump{\bV}^\mathrm{T}_{\xr}\tilde{\bm{F}}_{\xr}=\jump{\psi}_{\xr},
  \end{equation}
  with
$\jump{a}_{\xr}:=a_{i+1}-a_{i}$ and $\mean{a}_{\xr}:=\dfrac12(a_i+a_{i+1})$,
  then the scheme \eqref{eq:semi} with   $\tilde{\bF}_{\xr}$ is
  second-order accurate and entropy conservative.
  The corresponding numerical entropy flux is
  $\tilde{q}_{\xr}=\mean{\bV}^\mathrm{T}_{\xr}\tilde{\bF}_{\xr}-\mean{\psi}_{\xr}$.
\end{theorem}
For the scalar equation,  solving \eqref{eq:conserFlux} can uniquely give
  $\tilde{\bm{F}}_{\xr}$, but it is
not clear for a general system. In \cite{Tadmor1987},
 a solution of \eqref{eq:conserFlux} was constructed by the  path integral in
the phase space
\begin{equation}
  \tilde{\bF}_{\xr}=\int_{0}^{1} \bF(\bV_i+\xi(\bV_{i+1}-\bV_i))\dd \xi,
\end{equation}
which might be very hard to calculate except in some special
cases\cite{Fjordholm2009}. An explicit solution of \eqref{eq:conserFlux} was
given in \cite{Tadmor2004}, but it was
both expensive and numerically unstable. Some explicit algebraic
solutions of \eqref{eq:conserFlux}  were constructed in the literature for the specific systems,
such as the linear symmetric system, the shallow water equations \cite{Fjordholm2009},
the Euler equations \cite{Ismail2009}.
Herein, we  construct the {affordable} entropy conservative flux for the one-dimensional special RHD
equations using the strategy introduced in \cite{Ranocha}.
The key is to use the identity
\begin{equation*}
  \jump{ab}=\mean{a}\jump{b}+\mean{b}\jump{a},
\end{equation*}
where $\jump{a}$ and $\mean{a}$ denote the jump and mean of $a$, respectively,
and rewrite the jumps of the entropy variables ${\bV}$ and  the
potential ${\psi}$ as some {linear} combinations of the jump of a specially chosen
parameter vector.
To be specific, we first deal with the Lorentz factor and omit the subscript $i$  for simplicity.
Because
\begin{equation*}
  \dfrac{1}{\sqrt{1-u_R^2}}-\dfrac{1}{\sqrt{1-u_L^2}}
  =\dfrac{(u_L+u_R)(u_R-u_L)}{\sqrt{1-u_L^2}\sqrt{1-u_R^2}(\sqrt{1-u_L^2}+\sqrt{1-u_R^2})},
\end{equation*}
where the subscripts $L$ and $R$ denote the left and right values of the
variables used to calculate the entropy conservative flux,
one can define the ``Lorentz mean'' by
\begin{equation*}
  \jump{W}=\meanL{u}\jump{u},
\end{equation*}
where
\begin{equation}\label{eq:1Ddef}
  \meanL{u}
  =\dfrac{u_L+u_R}{\sqrt{1-u_L^2}\sqrt{1-u_R^2}(\sqrt{1-u_L^2}+\sqrt{1-u_R^2})}.
\end{equation}
If choosing the parameter vector
\begin{align*}
  &\bz=(z_1,z_2,z_3)^\mathrm{T}= (\rho,\rho/p,u)^\mathrm{T},
\end{align*}
then
\begin{equation}\label{eq:Vjump}
  \left\{
\begin{aligned}
  &\jump{\bV_1}=\dfrac{\jump{z_1}}{\meanln{z_1}}+\dfrac{1}{\Gamma-1}\dfrac{\jump{z_2}}{\meanln{z_2}}+\jump{z_2},\\
  &\jump{\bV_2}=\mean{uW}\jump{z_2}+\mean{z_2}\mean{W}\jump{z_3}+{\mean{z_2}\mean{z_3}{\meanL{z_3}}\jump{z_3}},\\
  &\jump{\bV_3}=-\mean{W}\jump{z_2}-\mean{z_2}{\meanL{z_3}}{\jump{z_3}},\\
  &\jump{\psi}=\mean{uW}\jump{z_1}+\mean{z_1}\mean{W}\jump{z_3}+\mean{z_1}\mean{z_3}{\meanL{z_3}}{\jump{z_3}},
\end{aligned}
\right.
\end{equation}
where $\meanln{a}= {\jump{a}}/{\jump{\ln{a}}}$ is the logarithmic mean
introduced in \cite{Ismail2009}, where its stable numerical implementation can be found.
If assuming that the entropy conservative flux is $\tilde{\bF}=(\tilde{\bF}_1,
\tilde{\bF}_2, \tilde{\bF}_3)^\mathrm{T}$,
and substituting the equations  \eqref{eq:Vjump} in \eqref{eq:conserFlux},
then
\begin{equation*}
  \left\{
\begin{aligned}
  &\dfrac{\tilde{\bF}_1}{\meanln{z_1}}=\mean{uW},\\
  &\dfrac{\tilde{\bF}_1}{(\Gamma-1)\meanln{z_2}}+\tilde{\bF}_1
  +\mean{uW}{\tilde{\bF}_2}-\mean{W}\tilde{\bF}_3=0,\\
  &\tilde{\bF}_2\mean{z_2}(\mean{W}
    + \mean{z_3}{\meanL{z_3}})
  - {{\tilde{\bF}_3}\mean{z_2}}{\meanL{z_3}}=
  \mean{z_1}(\mean{W}+{\mean{z_3}}{\meanL{z_3}}).
\end{aligned}
\right.
\end{equation*}
Solving the above equations  can obtain the entropy conservative flux for the one-dimensional special RHD
equations as follows
\begin{equation}\label{eq:ecFlux1D}
  \left\{
\begin{aligned}
  \tilde{\bF}_1=&\meanln{z_1}\mean{uW},\\
  \tilde{\bF}_2=&Q^{-1}\Bigg[ \left(1+\dfrac{1}{(\Gamma-1)\meanln{z_2}}\right){\mean{z_2}}{\meanL{z_3}}\tilde{\bF}_1
  + \mean{z_1}\mean{W}^2  \\
&+{\mean{z_1}\mean{z_3}\mean{W}}{\meanL{z_3}} \Bigg],\\
  \tilde{\bF}_3=&Q^{-1}\Bigg[ \mean{z_1}\mean{W}\mean{uW}+{\mean{z_1}\mean{z_3}\mean{uW}}{\meanL{z_3}} \\
    &+\left(1+\dfrac{1}{(\Gamma-1)\meanln{z_2}}\right)\tilde{\bF}_1
  \left(\mean{z_2}\mean{W}+{\mean{z_2}\mean{z_3}}{\meanL{z_3}} \right)\Bigg],
\end{aligned}
\right.
\end{equation}
where
$Q=\mean{z_2}\mean{W}^2+{\mean{z_2}\mean{z_3}\mean{W}}{\meanL{z_3}}
-{\mean{z_2}\mean{uW}}{\meanL{z_3}}$.
\begin{remark}
  It is worth noting that $Q$ in the above expressions is positive
 so that our entropy conservative flux is well defined. In fact,
 if $u_L=u_R$ and $W_L=W_R$, then
      $\meanL{z_3}=\meanL{u}=u_LW_L^3$ due to \eqref{eq:1Ddef}, thus $Q=\mean{z_2}W_L^2>0$;
 otherwise, $\meanL{z_3}=\meanL{u}= {\jump{W}}/{\jump{u}}$,
      thus
  \begin{align*}
    Q&=\mean{z_2}\left\{\dfrac{(W_L+W_R)^2}{4}
      +\left[\dfrac{u_L+u_R}{2}\dfrac{W_L+W_R}{2}-\dfrac{u_LW_L+u_RW_R}{2}\right]
      \dfrac{W_L-W_R}{u_L-u_R} \right\}\\
      &=\mean{z_2}W_LW_R>0.
  \end{align*}
\end{remark}

\begin{remark}
  It is  also easy to verify that the entropy conservative flux
  \eqref{eq:ecFlux1D} is consistent with the flux $\bF$.
  If letting $(\rho_L,u_L,p_L)=(\rho_R,u_R,p_R)=(\rho,u,p)$,
  then
  \begin{align*}
    \tilde{\bF}_1=&\rho {uW},\\
    \tilde{\bF}_2=&\dfrac{ \left(1+\dfrac{p}{(\Gamma-1)\rho}\right)\dfrac{{\rho^2}{u^2W^4}}{p}
  + \rho{W}^2 +{\rho u^2W^4} }{\rho W^2/p}
  =\rho hW^2 u^2+p,\\
    \tilde{\bF}_3=&\dfrac{\rho W^2u+\rho W^4u^3+
      \left(1+\dfrac{p}{(\Gamma-1)\rho}\right)\dfrac{{\rho^2 W^2u+\rho^2W^4 u^3}}{p} }{\rho W^2/p}
      =\rho hW^2u.
  \end{align*}
\end{remark}

The scheme \eqref{eq:semi} with the entropy conservative flux  \eqref{eq:ecFlux1D} is only
second-order accurate. However, if
using that entropy conservative flux as a building block, then one can obtain
an entropy conservative flux of the $2p$th-order ($p\in \bbN^+$) accurate scheme  by
using the linear combinations of the ``second-order accurate'' entropy conservative fluxes
\cite{Lefloch2002}. Here only presents the specific expressions for the ``$6$th-order accurate'' entropy conservative flux
\begin{align}
  \tilde{\bF}^{6\mbox{\scriptsize th}}_{\xr}=&\dfrac32\tilde{\bF}(\bU_i,\bU_{i+1})-\dfrac{3}{10}\left(\tilde{\bF}(\bU_{i-1},\bU_{i+1})+\tilde{\bF}(\bU_i,\bU_{i+2})\right)\nonumber\\
  &+\dfrac{1}{30}\left(\tilde{\bF}(\bU_{i-2},\bU_{i+1})+\tilde{\bF}(\bU_{i-1},\bU_{i+2})+\tilde{\bF}(\bU_i,\bU_{i+3})\right).
\end{align}
The readers are referred to  \cite{Lefloch2002,Fjordholm2012} for more details on constructing the ``high-order accurate'' entropy conservative flux.

\subsection{Entropy stable fluxes}\label{sec:2.2}
The entropy of hyperbolic conservation laws is conserved only if the solution is smooth.
In other words, the entropy is not conserved if the discontinuity such as the shock wave appears in the solution.
It is well-known that an entropy conservative scheme may become oscillatory near
the shock wave, thus we expect to construct an entropy stable scheme
by adding a dissipation term in the original entropy conservative
scheme. This section will first introduce the entropy stable flux
and its high-order extension of
Tadmor and his collaborators via the ENO reconstruction,
and then go to the low dissipative entropy stable
flux by using a switch function in the dissipation term \cite{Biswas2018}.
\begin{theorem}[Tadmor\cite{Tadmor1987}]\label{th:diff}
 If assuming that $\bm{D}_{\xr}$ is a symmetric positive semi-definite matrix
  and $\tilde{\bm{F}}_{\xr}$ is an entropy conservative flux,
  then the scheme \eqref{eq:semi} with the following numerical flux
  \begin{align}\label{eq:stableflux}
    \hat{\bF}_{\xr}=\tilde{\bF}_{\xr}-\dfrac12 \bm{D}_{\xr}\jump{\bV}_{\xr},
  \end{align}
  is entropy stable, i.e., satisfying the semi-discrete entropy inequality
  \begin{equation*}
    \dfrac{\dd}{\dd t}\eta(\bU_i(t))+\dfrac{1}{\Delta
    x}\left(\hat{q}_{\xr}(t)-\hat{q}_{\xl}(t)\right)\leq0,
  \end{equation*}
  for some numerical entropy flux function $\hat{q}_{\xr}$ consistent with $q$.
\end{theorem}
The above theorem holds for any positive semi-definite matrix $\bm{D}_{\xr}$, which
is usually chosen as
\begin{equation*}
  \bm{D}_{\xr}=\bm{R}_{\xr}\abs{\bm{\Lambda}_{\xr}}\bm{R}_{\xr}^\mathrm{T},
\end{equation*}
where $\bm{R}$ is a scaled  matrix of right eigenvectors,
whose existence can be ensured by the eigenvector scaling theorem in
\cite{Merriam1989}, and   satisfies
\begin{equation*}
\pd{\bF}{\bU}=\bm{R}\bm{\Lambda}\bm{R}^{-1},\quad
\bU_{\bV}=\bm{R}\bm{R}^\mathrm{T}.
\end{equation*}
For the one-dimensional special RHD equations, after some algebraic manipulations, the
  scaled matrix $\bm{R}$ is
\begin{align}\label{eq:scale1D}
  \begin{bmatrix}
    1 & 1 & 1 \\
    (u-c_s)Wh & uW & (u+c_s)Wh \\
    (1-uc_s)Wh & W & (1+uc_s)Wh \\
  \end{bmatrix}
  \begin{bmatrix}
    \frac{\rho W(1-u c_s)}{2\Gamma} & 0 & 0 \\
    0 & \frac{(\Gamma-1)\rho W}{\Gamma} & 0 \\
    0 & 0 & \frac{\rho W(1+u c_s)}{2\Gamma} \\
  \end{bmatrix}
  ^{\frac12},
\end{align}
where $c_s=\sqrt{{\Gamma p}/{\rho h}}$ is the sound speed.
Because the scaled matrix $\bm{R}$ is mainly defined at $x_{\xr}$, one has to use some ``averaged'' values
of the variables to calculate it. This paper chooses
\begin{equation*}
  \bar{\rho}=\meanln{\rho}_{\xr},\quad
  \bar{u}=\mean{u}_{\xr},\quad
  \bar{p}=\dfrac{\meanln{\rho}_{\xr}}{\meanln{\rho/p}_{\xr}},
\end{equation*}
to replace the variables $\rho,u,p$ in \eqref{eq:scale1D} to obtain the scaled
matrix $\bm{R}_{\xr}$.

For the choice of $\abs{\bm{\Lambda}}$,
one can use the Roe type dissipation term
\begin{equation*}
  \abs{\bm{\Lambda}}=\text{diag}\{\abs{\lambda_1},  \abs{\lambda_2}, \abs{\lambda_3}\},
\end{equation*}
where $\lambda_1,\lambda_2,\lambda_3$ are three eigenvalues of $\pd{\bF}{\bU}$,
or the Lax-Friedrichs type dissipation term
\begin{equation*}
  \abs{\bm{\Lambda}}=\max\{\abs{\lambda_1},\dots,\abs{\lambda_m}\}\bm{I}.
\end{equation*}

If only $\bV_i$ and $\bV_{i+1}$  are used to calculate $\jump{\bV}_{i+\frac12}$,
then the scheme \eqref{eq:semi} with the   entropy stable flux
\eqref{eq:stableflux} is only first-order accurate
 even if a ``high-order accurate'' entropy conservative flux is used.
In \cite{Fjordholm2012}, the arbitrary high-order entropy stable
schemes are constructed by  applying the ENO reconstruction to the scaled entropy variables
$\bw=\bm{R}^\mathrm{T}\bV$. More specifically, apply the $k$th-order accurate ENO reconstruction to
 $\bw$ to obtain the left and right limit values at $x_{\xr}$, denoted by
 $\bw_{\xr}^-$ and $\bw_{\xr}^+$, and
\begin{equation*}
  \jumpangle{\bw}_{\xr}=\bw_{\xr}^+-\bw_{\xr}^-,
\end{equation*}
  then replace the second-order entropy conservative flux \eqref{eq:ecFlux1D}
  with $2p$th-order entropy conservative flux $\tilde{\bF}^{2p\mbox{\scriptsize th}}$,
  where $p=k/2$ for even $k$ and $p=(k+1)/2$ for odd $k$,
  and replace the dissipation term with $\dfrac12
  \bm{R}_{\xr}\abs{\bm{\Lambda}_{\xr}}\jumpangle{\bw}_{\xr}$.
%
Finally, one has the ``$k$th-order accurate'' entropy stable flux
\begin{equation}\label{eq:HOstable}
  \hat{\bF}_{\xr}=\tilde{\bF}^{2p\mbox{\scriptsize th}}_{\xr}-\dfrac12
  \bm{R}_{\xr}\abs{\bm{\Lambda}_{\xr}}\jumpangle{\bw}_{\xr}.
\end{equation}
The semi-discrete numerical schemes \eqref{eq:semi} with
above high-order flux is entropy stable if the reconstruction satisfies the
following ``sign'' property \cite{Fjordholm2012}
\begin{equation*}
  \text{sign}(\jumpangle{\bw}_{\xr})=\text{sign}(\jump{\bw}_{\xr}),
\end{equation*}
which is satisfied by the ENO reconstructions \cite{Fjordholm2013}.

Certainly, one can also obtain higher-order accuracy with the WENO reconstruction instead of
the ENO reconstruction if  the same number of candidate points values are used,
but a general WENO reconstruction may not satisfy the ``sign'' property.
Borrowing the idea from \cite{Biswas2018}, we add a switch
function in the dissipation term as follows
\begin{equation}\label{eq:switch}
  \hat{\bm{F}}_{\xr}=\tilde{\bF}^{2pth}_{\xr}-\dfrac12
  \bm{S}_{\xr}\bm{R}_{\xr}\abs{\bm{\Lambda}_{\xr}}\jumpangle{\bw}_{\xr},
\end{equation}
where
\begin{equation*}
  \bm{S}^l_{\xr}=\begin{cases}
    1,\quad &\text{if}
    ~\text{sign}(\jumpangle{\bw}^l_{\xr})=\text{sign}(\jump{\bw}^l_{\xr})\not= 0, \\
    0,\quad &\text{otherwise},
  \end{cases}
\end{equation*}
here the superscript $l$ denotes the $l$-th entry of the diagonal matrix $\bm{S}_{\xr}$
or the $l$-th component of the jump of $\bw$.
When the WENO reconstruction does not satisfy the ``sign'' property,
corresponding dissipation term becomes zero, and thus the semi-discrete numerical
scheme with the flux \eqref{eq:switch} is entropy stable according to Theorem
\ref{th:diff}.
Meanwhile, compared to the entropy stable flux using the ENO reconstruction, the
above flux using the WENO reconstruction leads to less dissipation because the switch function is
not active  at all locations.

This paper  uses the fifth-order accurate WENO reconstruction in \cite{jiang}
and
the following  third-order accurate  Runge-Kutta time discretization for
the time derivatives in \eqref{eq:semi}
\begin{align*}
  &\bU^{(1)}=\bU^n+\Delta t \bm{L}(\bU^n),\\
  &\bU^{(2)}=\dfrac34\bU^n+\dfrac14\left(\bU^{(1)}+\Delta t \bm{L}(\bU^{(1)})\right),\\
  &\bU^{n+1}=\dfrac13\bU^n+\dfrac23\left(\bU^{(2)}+\Delta t \bm{L}(\bU^{(2)})\right),
\end{align*}
where $[\bm{L}(\bU)]_i$ denotes the right-hand side term of \eqref{eq:semi}.

\section{Two-dimensional  schemes}\label{section:2D}
The two-dimensional finite difference scheme for solving \eqref{eq:RHD} with
$d=2$ can be done in a dimension-by-dimension fashion.
For simplicity, the notations $\bF_1,\bF_2,u_1,u_2,m_1,m_2$ and $x_1,x_2$ are
replaced with $\bF,\bG,u,v,m_x,m_y$ and $x,y$ respectively,
and thus $\bF=(Du,m_xu+p,m_yu,m_x)^\mathrm{T}, \bG=(Dv,m_xv,m_yv+p,m_y)^\mathrm{T}$.
The two-dimensional entropy pair and potential are
$$\eta=\dfrac{-\rho WS}{\Gamma-1},\quad q_x=\dfrac{-\rho u WS}{\Gamma-1},\quad
q_y=\dfrac{-\rho v WS}{\Gamma-1},\quad \psi_x=\rho uW, \quad \psi_y=\rho vW.$$
Consider a uniform Cartesian mesh with the spatial stepsizes $\Delta x,\Delta y$.
The solution $\bU$ is approximated at $(x_i,y_j),i=1,\cdots,N_x,j=1,\cdots,N_y$, and the
$x$- and $y$-directional numerical fluxes are defined at
$(x_{\xr},y_j)$ and $(x_i,y_{\yr})$, respectively.
Then a semi-discrete conservative finite
difference scheme can be expressed as
\begin{equation}\label{eq:semi2D}
  \dfrac{\dd}{\dd t}\bU_{i,j}(t)=
  -\dfrac{1}{\Delta x}\left(\hat{\bF}_{\xr,j}-\hat{\bF}_{\xl,j}\right)
  -\dfrac{1}{\Delta y}\left(\hat{\bG}_{i,\yr}-\hat{\bG}_{i,\yl}\right),
\end{equation}
where the numerical fluxes $\hat{\bF}_{\xr,j},\hat{\bG}_{i,\yr}$ are  defined by
\begin{align*}
  &\hat{\bm{F}}_{\xr,j}=\tilde{\bF}^{2p\mbox{\scriptsize th}}_{\xr,j}-\dfrac12
  \bm{S}_{\xr,j}\bm{R}^x_{\xr,j}\abs{\bm{\Lambda}^x_{\xr,j}}\jumpangle{\bw^x}_{\xr,j},\\
  &\hat{\bm{G}}_{i,\yr}=\tilde{\bG}^{2p\mbox{\scriptsize th}}_{i,\yr}-\dfrac12
  \bm{S}_{i,\yr}\bm{R}^y_{i,\yr}\abs{\bm{\Lambda}^y_{i,\yr}}\jumpangle{\bw^y}_{i,\yr},
\end{align*}
with $\bw^x=\bm{R}^x \bV, \bw^y=\bm{R}^y \bV$.
Here $\tilde{\bF}^{2p\mbox{\scriptsize th}},\tilde{\bG}^{2p\mbox{\scriptsize th}}$, $\bm{R}^x,\bm{R}^y$,
and $\bm{\Lambda}^x,\bm{\Lambda}^y$ are the ``high-order accurate'' entropy
conservative fluxes,
the scaled matrices of right eigenvectors, and the diagonal Roe type or
Lax-Friedrichs type dissipation terms in $x$- and $y$-directions, respectively,
which will be given below,
and $\bm{S}$ is the same switch function as in the one-dimensional case.

Motivated by the one-dimensional case, one can define two ``Lorentz mean'' by
\begin{equation*}
  \jump{W}=\meanLx{u,v}\jump{u}+\meanLy{u,v}\jump{v},
\end{equation*}
where
\begin{align*}
  &\meanLx{u,v}=\dfrac{u_L+u_R}{\sqrt{1-u_L^2-v_L^2}\sqrt{1-u_R^2-v_R^2}
  (\sqrt{1-u_L^2-v_L^2}+\sqrt{1-u_R^2-v_R^2})},\\
  &\meanLy{u,v}=\dfrac{v_L+v_R}{\sqrt{1-u_L^2-v_L^2}\sqrt{1-u_R^2-v_R^2}
  (\sqrt{1-u_L^2-v_L^2}+\sqrt{1-u_R^2-v_R^2})}.
\end{align*}
If taking the parameter vector
$  \bz=(z_1,z_2,z_3,z_4)^\mathrm{T}=(\rho,\rho/p,u,v)^\mathrm{T}$
and following the same procedure in the one-dimensional case,
then one can obtain the entropy conservative flux $\tilde{\bF}=(\tilde{\bF}_1,
\tilde{\bF}_2, \tilde{\bF}_3, \tilde{\bF}_4)^\mathrm{T}$ in the $x$-direction
\begin{equation*}
  \left\{
\begin{aligned}
  \tilde{\bF}_1=&\meanln{z_1}\mean{uW},\\
  \tilde{\bF}_2=&Q^{-1}\Big\{ \alpha\mean{z_2}\meanLx{z_3,z_4}\tilde{\bF}_1
    +\mean{z_1}(\mean{W}^2-\mean{vW}\meanLy{z_3,z_4}) \\
  &+\mean{z_1}\mean{W}(\mean{z_3}\meanLx{z_3,z_4}+\mean{z_4}\meanLy{z_3,z_4}) \Big\},\\
  \tilde{\bF}_3=&Q^{-1}\Big\{ \alpha\mean{z_2}\meanLy{z_3,z_4}\tilde{\bF}_1
    +\mean{z_1}\mean{uW}\meanLy{z_3,z_4}\Big\},\\
  \tilde{\bF}_4=&\mean{W}^{-1}\Big(\alpha\tilde{\bF}_1+\mean{uW}\tilde{\bF}_2+\mean{vW}\tilde{\bF}_3\Big),
\end{aligned}
\right.
\end{equation*}
and the entropy conservative flux $\tilde{\bG}=(\tilde{\bG}_1,
\tilde{\bG}_2, \tilde{\bG}_3, \tilde{\bG}_4)^\mathrm{T}$
in the $y$-direction
\begin{equation*}
  \left\{
\begin{aligned}
  \tilde{\bG}_1=&\meanln{z_1}\mean{vW},\\
  \tilde{\bG}_2=&Q^{-1}\Big\{ \alpha\mean{z_2}\meanLx{z_3,z_4}\tilde{\bG}_1
    +\mean{z_1}\mean{vW}\meanLx{z_3,z_4}\Big\},\\
    \tilde{\bG}_3=&Q^{-1}\Big\{ \alpha\mean{z_2}\meanLy{z_3,z_4}\tilde{\bG}_1
      +\mean{z_1}(\mean{W}^2-\mean{uW}\meanLx{z_3,z_4})\nonumber\\
    &+\mean{z_1}\mean{W}(\mean{z_3}\meanLx{z_3,z_4}+\mean{z_4}\meanLy{z_3,z_4}) \Big\},\\
  \tilde{\bG}_4=&\mean{W}^{-1}\Big(\alpha\tilde{\bG}_1+\mean{uW}\tilde{\bG}_2+\mean{vW}\tilde{\bG}_3\Big),
\end{aligned}
\right.
\end{equation*}
where
\begin{align}
  \alpha=&1+\dfrac{1}{(\Gamma-1)\meanln{z_2}},
  \\
  Q=&\mean{z_2}\mean{W}^2+\mean{z_2}\Big(\mean{z_3}\mean{W}\meanLx{z_3,z_4}
-\mean{uW}\meanLx{z_3,z_4}\nonumber\\
&+\mean{z_4}\mean{W}\meanLy{z_3,z_4}
-\mean{vW}\meanLy{z_3,z_4}\Big)\label{eq:Q2D}.
\end{align}
Similarly,  $Q$ in \eqref{eq:Q2D} is positive, because
if $W_L=W_R$, then
$\meanLx{z_3,z_4}=\dfrac{u_L+u_R}{2}W_L^3$, $\meanLy{z_3,z_4}=\dfrac{v_L+v_R}{2}W_L^3$,
and thus $Q=\mean{z_2}W_L^2>0$;
otherwise, $\meanLx{z_3,z_4}=\dfrac{(u_L+u_R)\jump{W}}{\jump{u^2}+\jump{v^2}}$,
$\meanLx{z_3,z_4}=\dfrac{(v_L+v_R)\jump{W}}{\jump{u^2}+\jump{v^2}}$,
and thus one can simplify $Q$ as $\mean{z_2}W_LW_R$, which is positive.
The above two-dimensional entropy conservative fluxes are also consistent
after some algebraic simplification.

For the two-dimensional special RHD equations, the scaled matrix $\bm{R}^x$ in $x$-direction is
\begin{align*}
  \begin{bmatrix}
    1                            & 1/W & Wv           & 1                             \\
    hW\mathcal{A}^x_{-}\lambda^x_{-} & u   & 2hW^2uv      & hW\mathcal{A}^x_{+}\lambda^x_{+}  \\
    hWv                          & v   & h(1+2W^2v^2) & hWv                           \\
    hW\mathcal{A}^x_{-}            & 1   & 2hW^2v       & hW\mathcal{A}^x_{+}             \\
  \end{bmatrix}
  \begin{bmatrix}
  \frac{\mathcal{B}^x-\mathcal{C}^x}{2} & 0 & 0 & 0 \\
  0 & \frac{(\Gamma-1)\rho W^3}{\Gamma} & 0 & 0 \\
  0 & 0 & \frac{p}{W(1-u^2)h} & 0 \\
  0 & 0 & 0 & \frac{\mathcal{B}^x+\mathcal{C}^x}{2} \\
  \end{bmatrix}
  ^{\frac12},
\end{align*}
where $\lambda^x_{\pm}=\dfrac{u(1-c_s^2)\pm
c_s/W\sqrt{1-u^2-v^2c_s^2}}{1-(u^2+v^2)c_s^2}$ are two eigenvalues in the $x$-direction,
and $\mathcal{A}^x_{\pm}=\dfrac{1-u^2}{1-u\lambda^x_{\pm}}$,
$\mathcal{B}^x=\dfrac{\rho W(1-u^2-v^2c_s^2)}{\Gamma(1-u^2)}$,
$\mathcal{C}^x=\dfrac{\rho uc_s\sqrt{1-u^2-v^2c_s^2}}{\Gamma(1-u^2)}$.
The scaled matrix $\bm{R}^y$ in $y$-direction is
\begin{align*}
  \begin{bmatrix}
    1                            & Wu           & 1/W & 1                             \\
    hWu                          & h(1+2W^2u^2) & u   & hWu                           \\
    hW\mathcal{A}^y_{-}\lambda^y_{-} & 2hW^2uv      & v   & hW\mathcal{A}^y_{+}\lambda^y_{+}  \\
    hW\mathcal{A}^y_{-}            & 2hW^2u       & 1   & hW\mathcal{A}^y_{+}             \\
  \end{bmatrix}
  \begin{bmatrix}
    \frac{\mathcal{B}^y-\mathcal{C}^y}{2} & 0 & 0 & 0 \\
    0 & \frac{p}{W(1-v^2)h} & 0 & 0 \\
    0 & 0 & \frac{(\Gamma-1)\rho W^3}{\Gamma} & 0 \\
    0 & 0 & 0 & \frac{\mathcal{B}^y+\mathcal{C}^y}{2} \\
  \end{bmatrix}
  ^{\frac12},
\end{align*}
where $\lambda^y_{\pm}=\dfrac{v(1-c_s^2)\pm
c_s/W\sqrt{1-v^2-u^2c_s^2}}{1-(u^2+v^2)c_s^2}$ are two eigenvalues in the $y$-direction,
and $\mathcal{A}^y_{\pm}=\dfrac{1-v^2}{1-v\lambda^y_{\pm}}$,
$\mathcal{B}^y=\dfrac{\rho W(1-v^2-u^2c_s^2)}{\Gamma(1-v^2)}$,
$\mathcal{C}^y=\dfrac{\rho vc_s\sqrt{1-v^2-u^2c_s^2}}{\Gamma(1-v^2)}$.

The dissipation terms $\bm{\Lambda}^x$ and $\bm{\Lambda}^y$ are similar to the one-dimensional
case {except for that} the eigenvalues used in the dissipation terms are
$\lambda^x_{-},u,u,\lambda^x_{+}$ for $\bm{\Lambda}^x$, and
$\lambda^y_{-},v,v,\lambda^y_{+}$ for $\bm{\Lambda}^y$, respectively.
In order to obtain the jumps $\jumpangle{\bw^x}$ and $\jumpangle{\bw^y}$,
one just needs to perform the WENO reconstructions in $x$- and $y$-directions
independently.
Similar to the one-dimensional case, the third-order Runge-Kutta scheme is
also used for the time derivatives in \eqref{eq:semi2D}.
This completes the description of the two-dimensional entropy stable finite
difference scheme for the special RHD equations.

\section{Numerical results}\label{section:Num}
This section presents some numerical results to
validate the performance of our entropy stable schemes for
 the special RHD equations \eqref{eq:RHD} with $d=1,2$.
All the tests take the CFL number as $0.4$, $\Gamma=5/3$, and the Lax-Friedrichs
type dissipation terms unless otherwise stated.
For the one- and two-dimensional tests,
the time stepsizes are respectively chosen as
\begin{equation*}
  \Delta t=
  \dfrac{\text{CFL}\Delta x}{\max\limits_{i}{\abs{\lambda^x(\bU_i)}}},
\end{equation*}
and
\begin{equation*}
  \Delta t=
  \dfrac{\text{CFL}}{\max\limits_{i,j}{\abs{\lambda^x(\bU_{i,j})}}/\Delta
  x+\max\limits_{i,j}{\abs{\lambda^y(\bU_{i,j})}}/\Delta y},
\end{equation*}
where $\abs{\lambda^x}$ and $\abs{\lambda^y}$ are the maximum absolute values of
all the eigenvalues in $x$- and $y$-directions, respectively, but
for the accuracy tests, $\Delta t$ is taken as the minimum between
the above choices and $\text{CFL}{\Delta x}^{5/3}$ (resp. $\text{CFL}\min(\Delta x,\Delta y)^{5/3}$)
for the one-dimensional (resp. two-dimensional) test to make the spatial error
{dominate}.

\subsection{One-dimensional results}
\begin{example}[Accuracy test]\label{ex:acc1D}\rm
	This test is used to verify the accuracy.
  The initial condition is
  \begin{equation*}
    (\rho,u,p)=(1+0.2\sin x,~0.2,~1), \quad x\in[0,2\pi],
  \end{equation*}
  with the periodic boundary condition.
  The exact solutions can be given by
  \begin{equation*}
    (\rho,u,p)=(1+0.2\sin (x-0.2t),~0.2,~1).
  \end{equation*}
\end{example}
Table \ref{tab:acc1D} lists the errors and the orders of convergence in $\rho$ at $t=0.1$
obtained by using our 1D scheme. It is seen that our scheme gets the fifth-order
accuracy as expected.

\begin{table}[!ht]
  \centering
  \begin{tabular}{r|cc|cc|cc} \hline
 $N$& $\ell^1$ error & order & $\ell^2$ error & order & $\ell^\infty$ error & order \\ \hline
 20 & 5.475e-06 &  -   & 6.741e-06 &  -   & 1.453e-05 &  -   \\
 40 & 1.615e-07 & 5.08 & 1.966e-07 & 5.10 & 3.979e-07 & 5.19 \\
 80 & 2.692e-09 & 5.91 & 3.450e-09 & 5.83 & 7.490e-09 & 5.73 \\
160 & 7.791e-11 & 5.11 & 1.054e-10 & 5.03 & 2.622e-10 & 4.84 \\
320 & 2.448e-12 & 4.99 & 3.331e-12 & 4.98 & 8.297e-12 & 4.98 \\
\hline
  \end{tabular}
  \caption{Example \ref{ex:acc1D}: Errors and orders of convergence of in $\rho$ at $t=0.1$.}
  \label{tab:acc1D}
\end{table}

\begin{figure}[!ht]
  \centering
  \subfigure[$\rho$]{
    \includegraphics[width=0.3\textwidth, trim=40 40 50 50, clip]{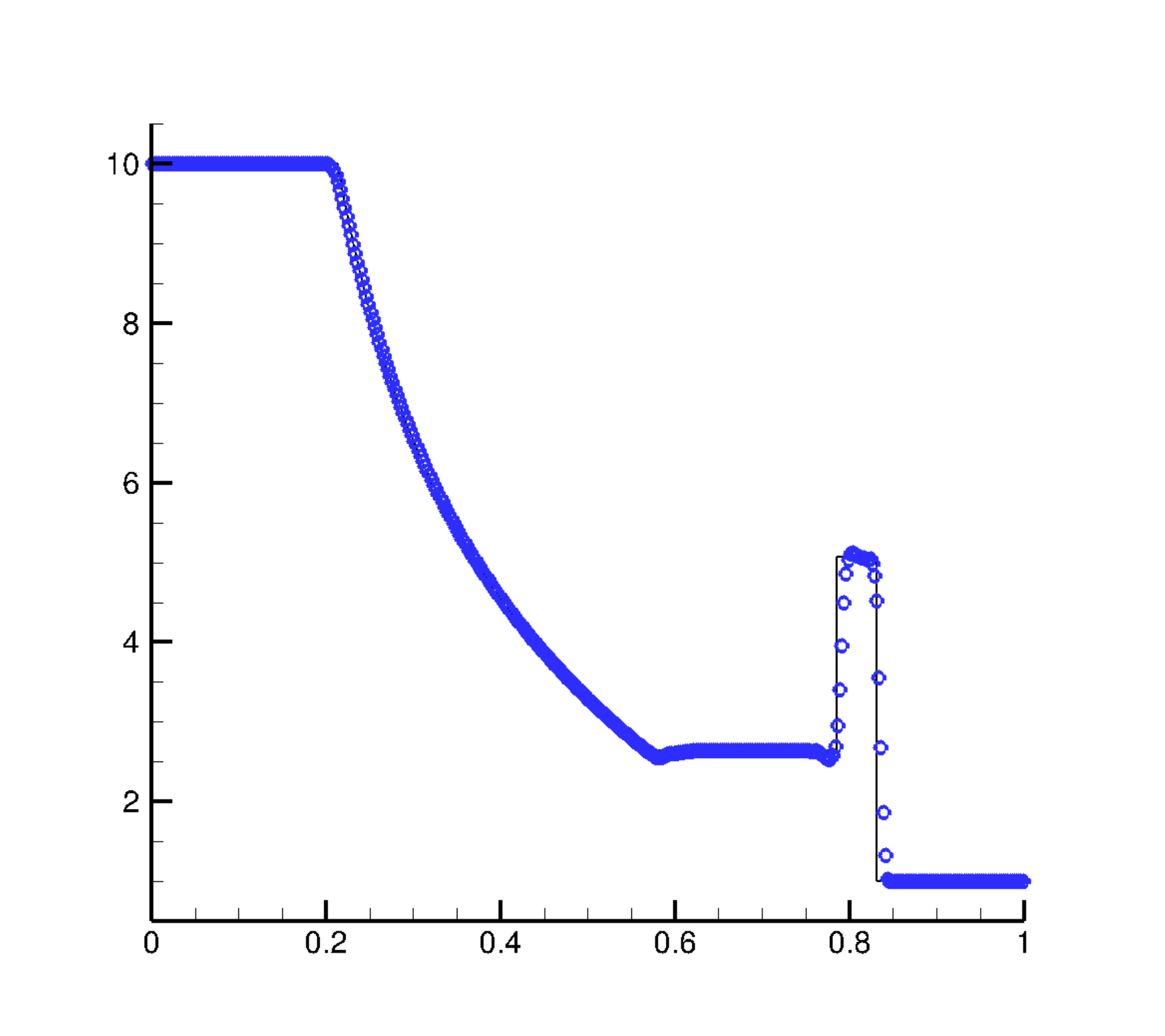}
  }
  \subfigure[$u$]{
    \includegraphics[width=0.3\textwidth, trim=40 40 50 50, clip]{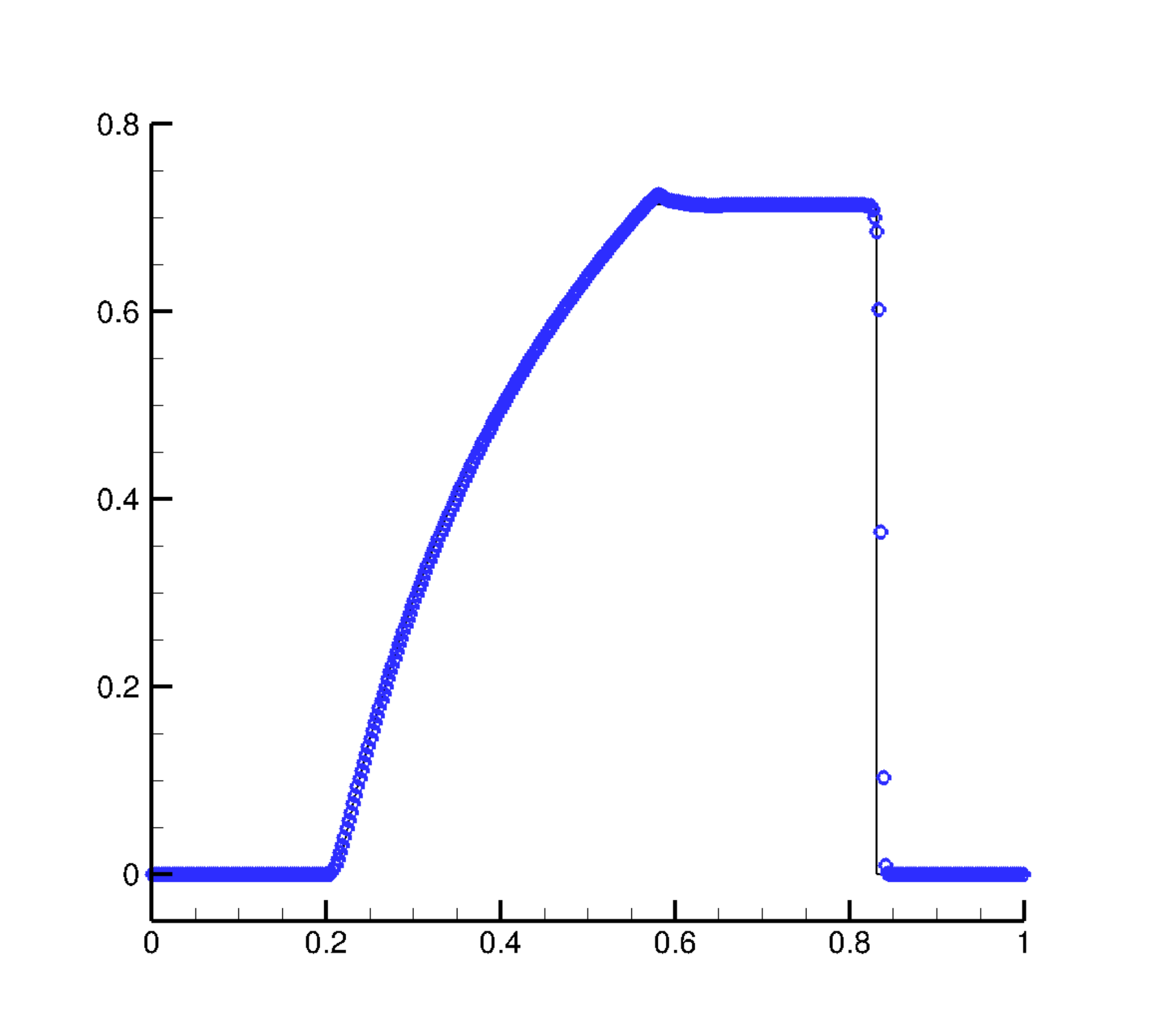}
  }
  \subfigure[$p$]{
    \includegraphics[width=0.3\textwidth, trim=40 40 50 50, clip]{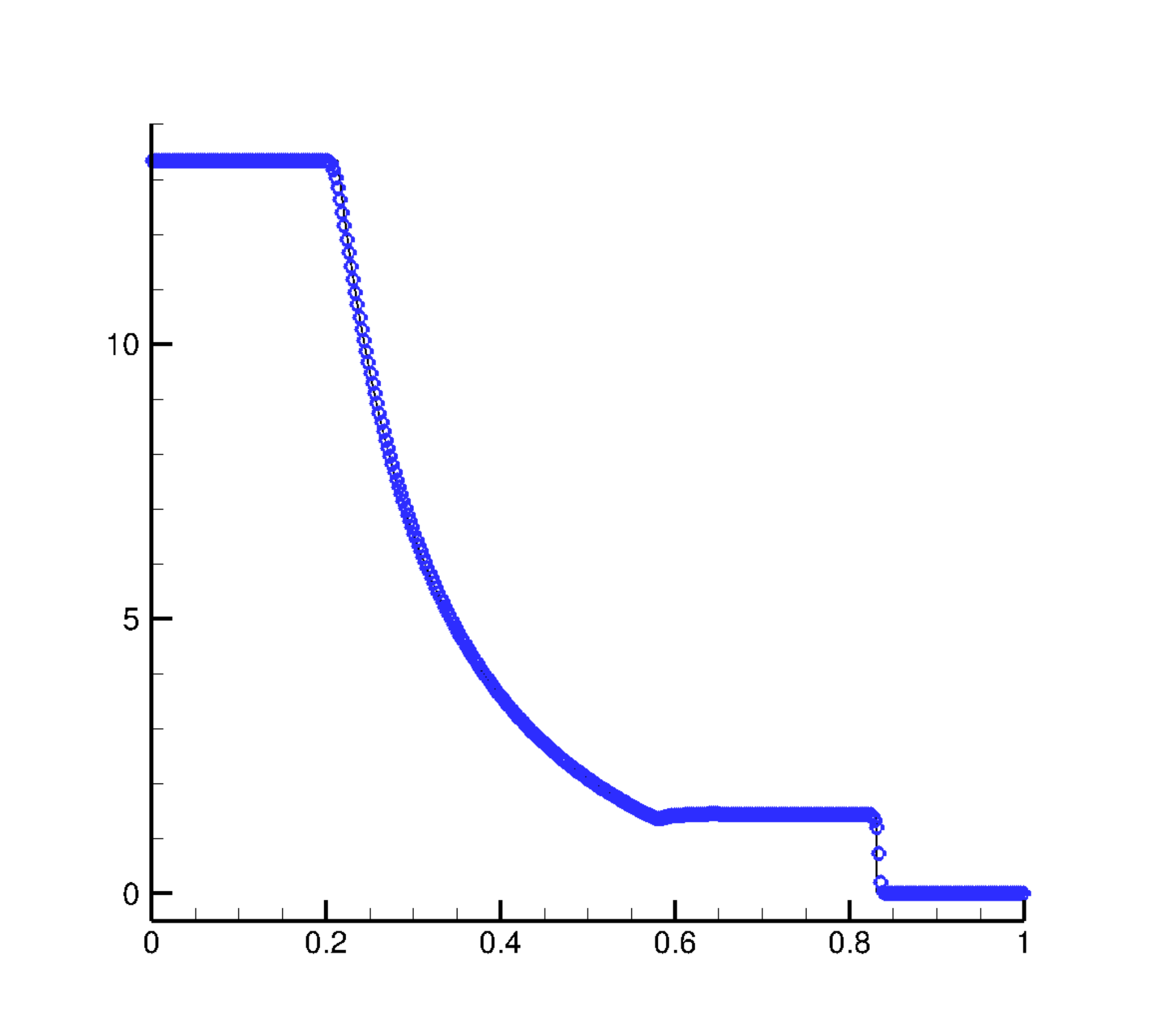}
  }
  \caption{Example \ref{ex:RP1}: Circles and
  the solid lines denote the numerical and  exact solutions, respectively. $N=400,~t=0.4$.}
  \label{fig:RP1}
\end{figure}

\begin{example}[Riemann problem 1]\label{ex:RP1}\rm
  The initial data of the first 1D Riemann problem are
  \begin{equation*}
    (\rho,u,p)=\begin{cases}
      (10,~0,~40/3),    &\quad x<0.5, \\
      (1,~0,~10^{-6}), &\quad x>0.5.
    \end{cases}
  \end{equation*}
\end{example}
As the time increases, the initial discontinuity will be decomposed into a
left-moving rarefaction wave, a contact discontinuity, and a right-moving shock
wave.
The rest-mass density, the velocity, and the pressure at $t=0.4$ obtained by the
entropy stable schemes with $400$ are shown in Figure \ref{fig:RP1}.
One can see that the numerical solutions are in good agreement with the exact
solutions, and the shock, the rarefaction wave, and the contact discontinuity are well
captured without obvious oscillations.

\begin{figure}[!ht]
  \centering
  \subfigure[$\rho$]{
    \includegraphics[width=0.3\textwidth, trim=40 40 50 50, clip]{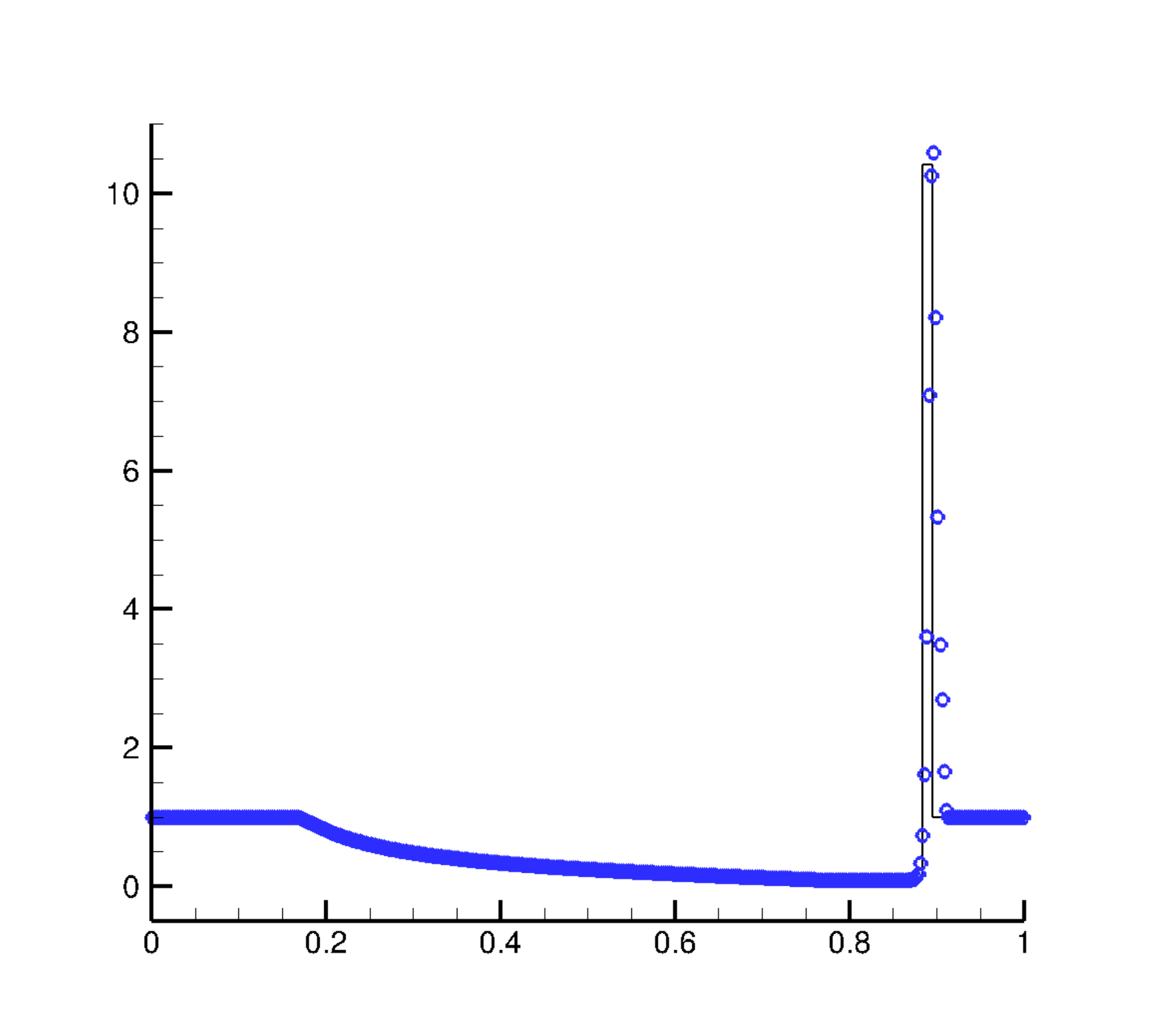}
  }
  \subfigure[$u$]{
    \includegraphics[width=0.3\textwidth, trim=40 40 50 50, clip]{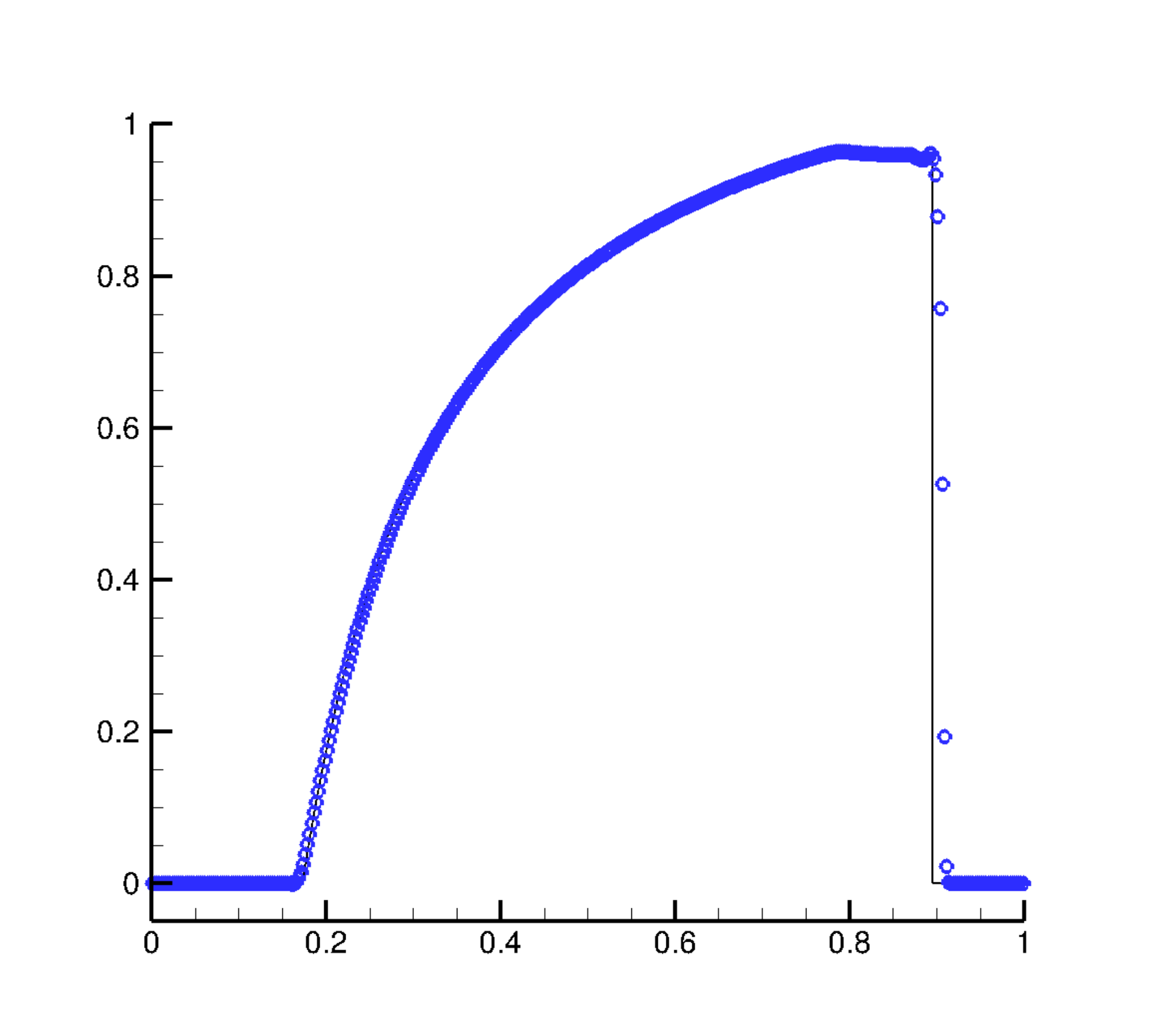}
  }
  \subfigure[$p$]{
    \includegraphics[width=0.3\textwidth, trim=40 40 50 50, clip]{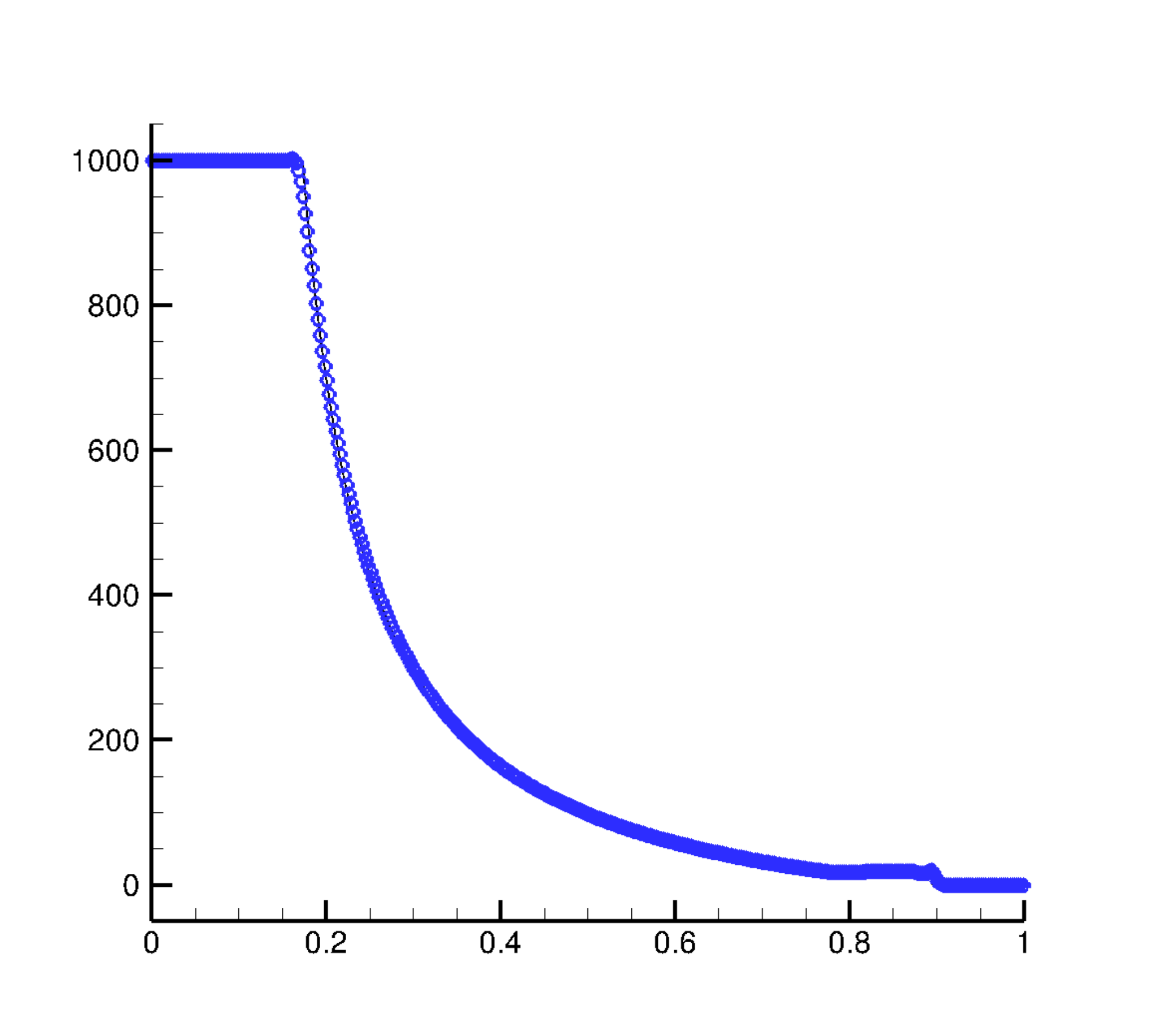}
  }
  \caption{Example \ref{ex:RP2}: Circles and
  the solid lines denote the numerical and  exact solutions, respectively. $N=400,~t=0.4$.}
  \label{fig:RP2}
\end{figure}

\begin{example}[Riemann problem 2]\label{ex:RP2}\rm
  The initial data of the second 1D Riemann problem are
  \begin{equation}
    (\rho,u,p)=\begin{cases}
      (1,~0,~10^3),    &\quad x<0.5, \\
      (1,~0,~10^{-2}), &\quad x>0.5.
    \end{cases}
  \end{equation}
\end{example}
The flow pattern is similar to that of the first Riemann problem, but more
extreme and difficult with a heavily curved profile for the rarefaction fan.
The region between the contact discontinuity and the right-moving
shock wave is very narrow so that resolving those waves is very challenging.
The rest-mass density, the velocity, and the pressure at $t=0.4$ obtained by the
entropy stable scheme with $400$ are shown in Figure \ref{fig:RP2}.
It can be seen that our scheme can still resolve the waves well, especially for the
rest-mass density profile, even though small undershoot appears in the narrow region
between the contact discontinuity and the right-moving
shock wave.

\begin{figure}[!ht]
  \centering
  \subfigure[$\rho$]{
    \includegraphics[width=0.3\textwidth, trim=40 40 50 50, clip]{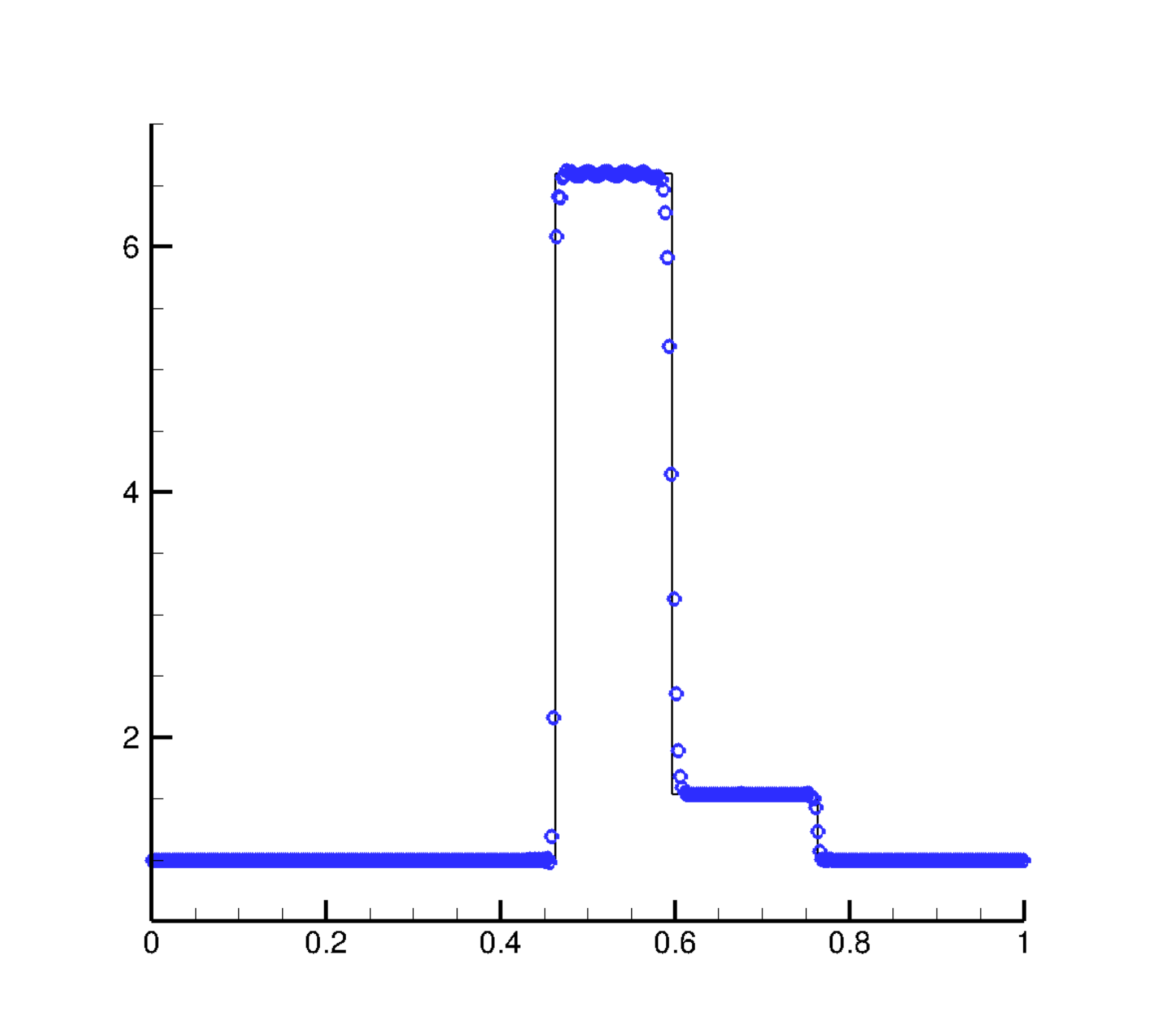}
  }
  \subfigure[$u$]{
    \includegraphics[width=0.3\textwidth, trim=40 40 50 50, clip]{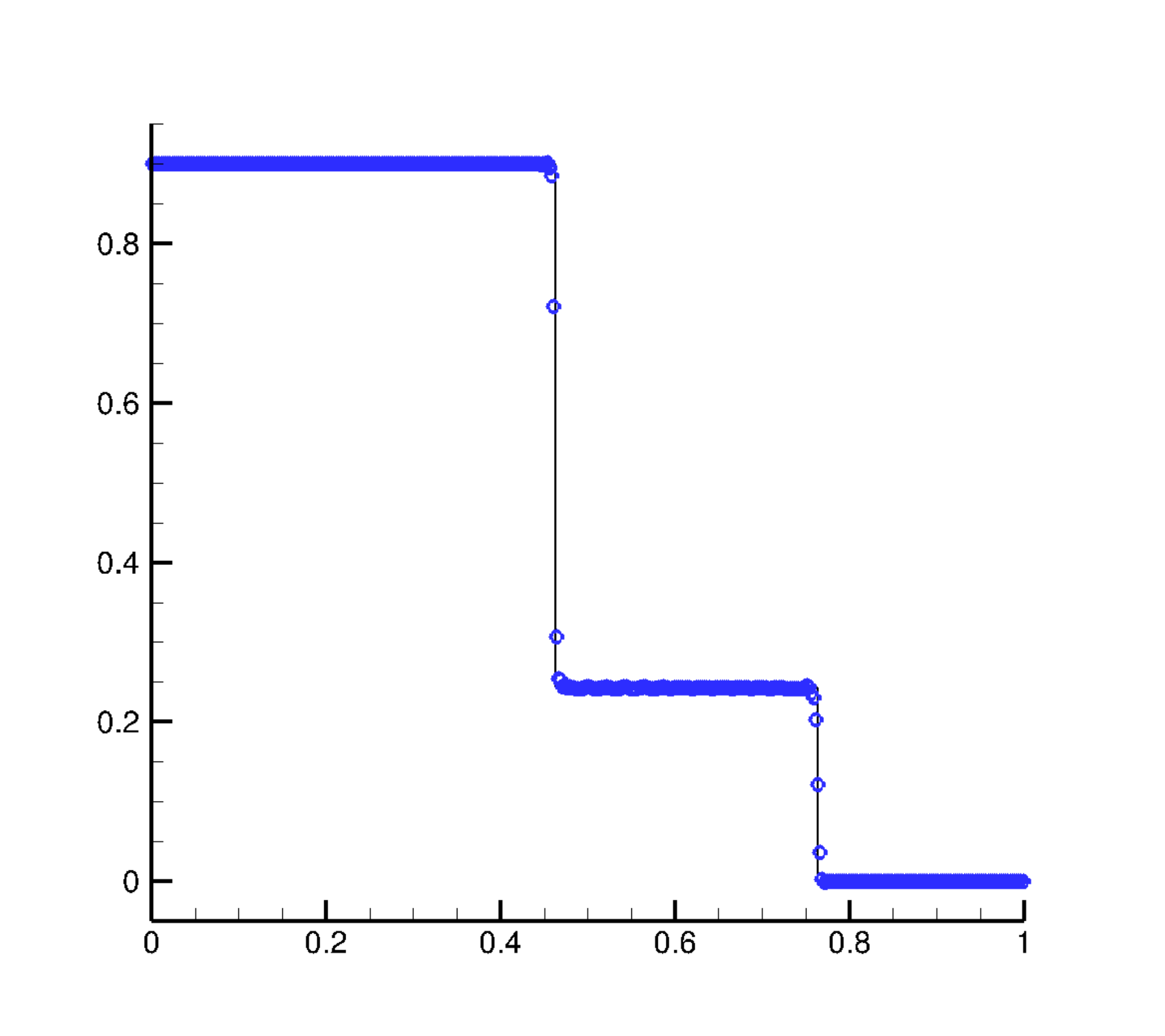}
  }
  \subfigure[$p$]{
    \includegraphics[width=0.3\textwidth, trim=40 40 50 50, clip]{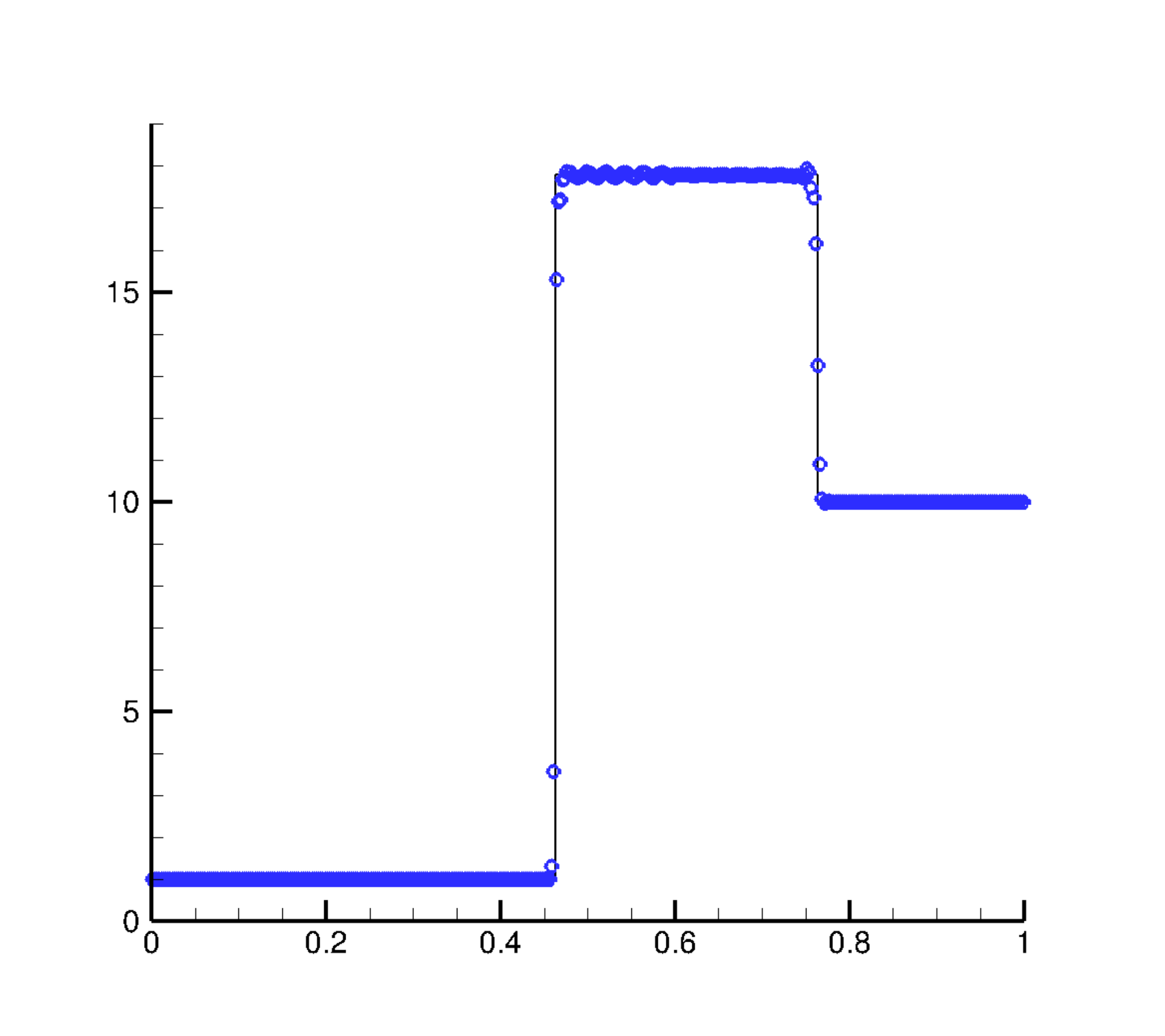}
  }
  \caption{Example \ref{ex:RP3}: Circles and
  the solid lines denote the numerical and  exact solutions, respectively. $N=400,~t=0.4$.}
  \label{fig:RP3}
\end{figure}

\begin{example}[Riemann problem 3]\label{ex:RP3}\rm
  The initial data of the third 1D Riemann problem are
  \begin{equation}
    (\rho,u,p)=\begin{cases}
      (1,~0.9,~1),    &\quad x<0.5, \\
      (1,~0,~10), &\quad x>0.5,
    \end{cases}
  \end{equation}
  with $\Gamma=4/3$.
\end{example}
The solutions will contain a slowly left-moving shock wave, a contact
discontinuity, and a fast right-moving shock wave, see Figure \ref{fig:RP3}, where
numerical solutions at $t=0.4$ are obtained by our entropy stable scheme with
$400$ uniform cells. Our numerical solutions are
in agreement with the exact solutions,
although small oscillations are observed behind the left-moving shock wave like
many shock-capturing methods,
but they can be
improved by using the adaptive moving mesh method, see \cite{he1}.

\begin{figure}[!ht]
  \centering
  \subfigure[$\rho$]{
    \includegraphics[width=0.3\textwidth, trim=40 40 50 50, clip]{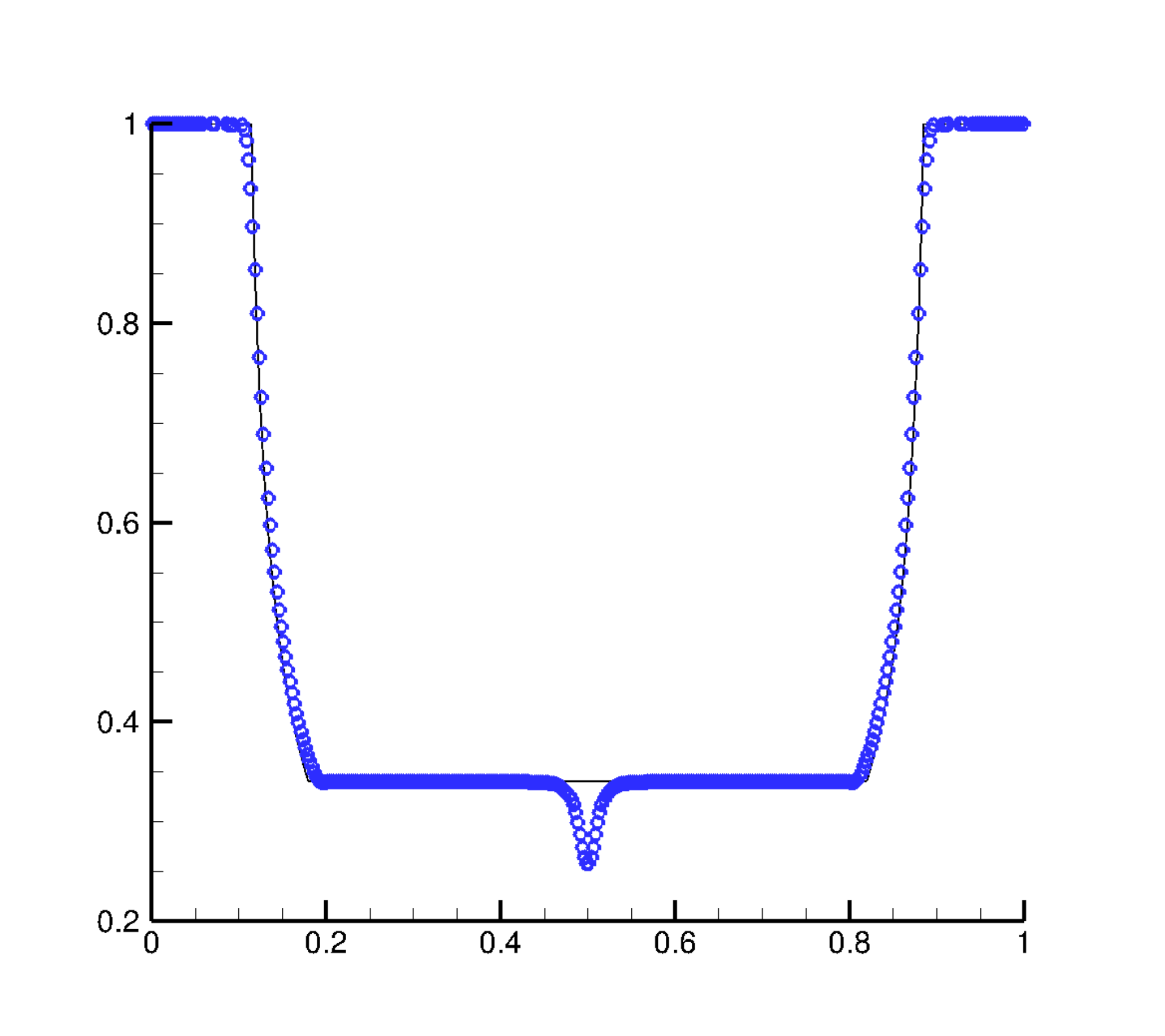}
  }
  \subfigure[$u$]{
    \includegraphics[width=0.3\textwidth, trim=40 40 50 50, clip]{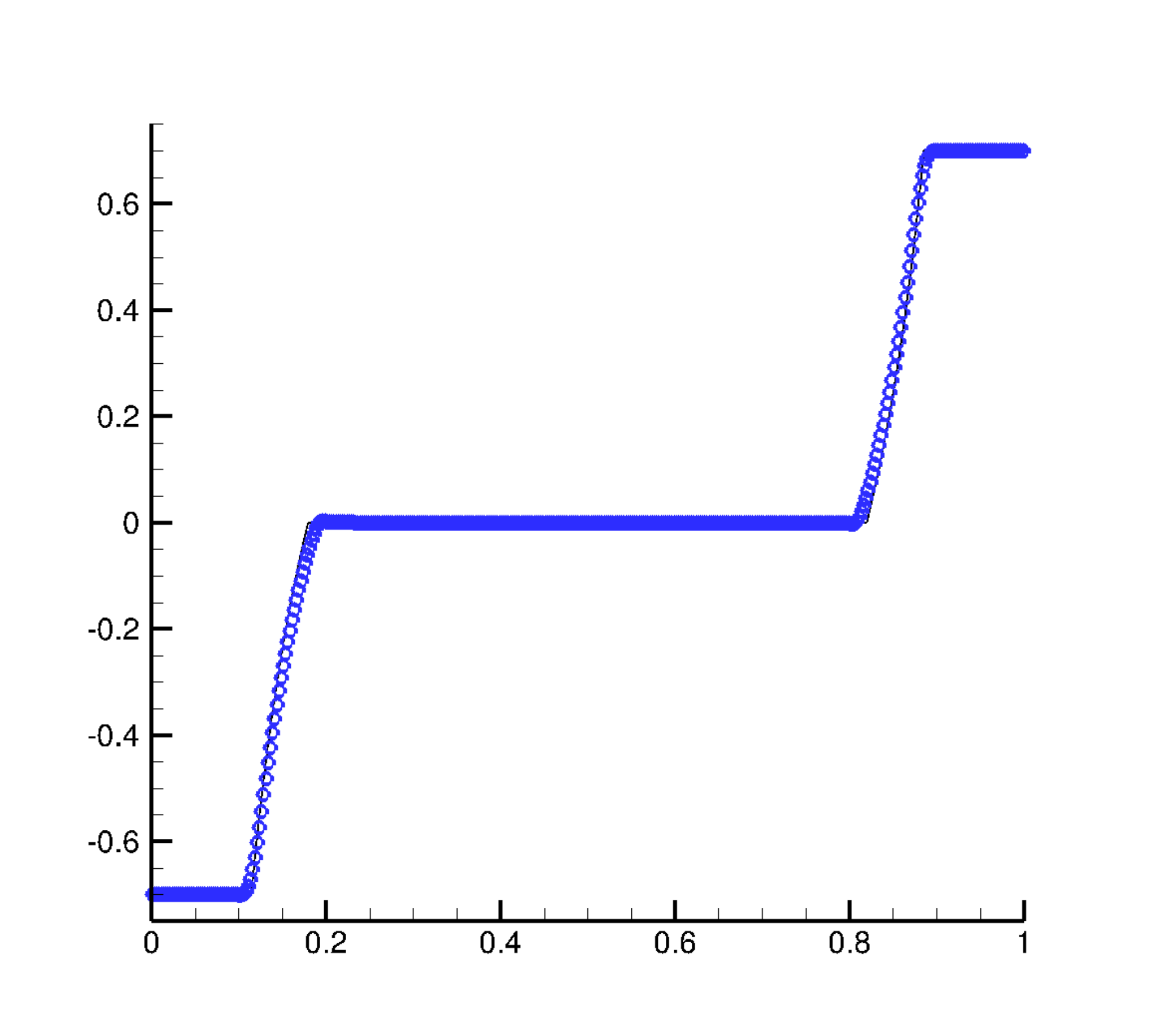}
  }
  \subfigure[$p$]{
    \includegraphics[width=0.3\textwidth, trim=40 40 50 50, clip]{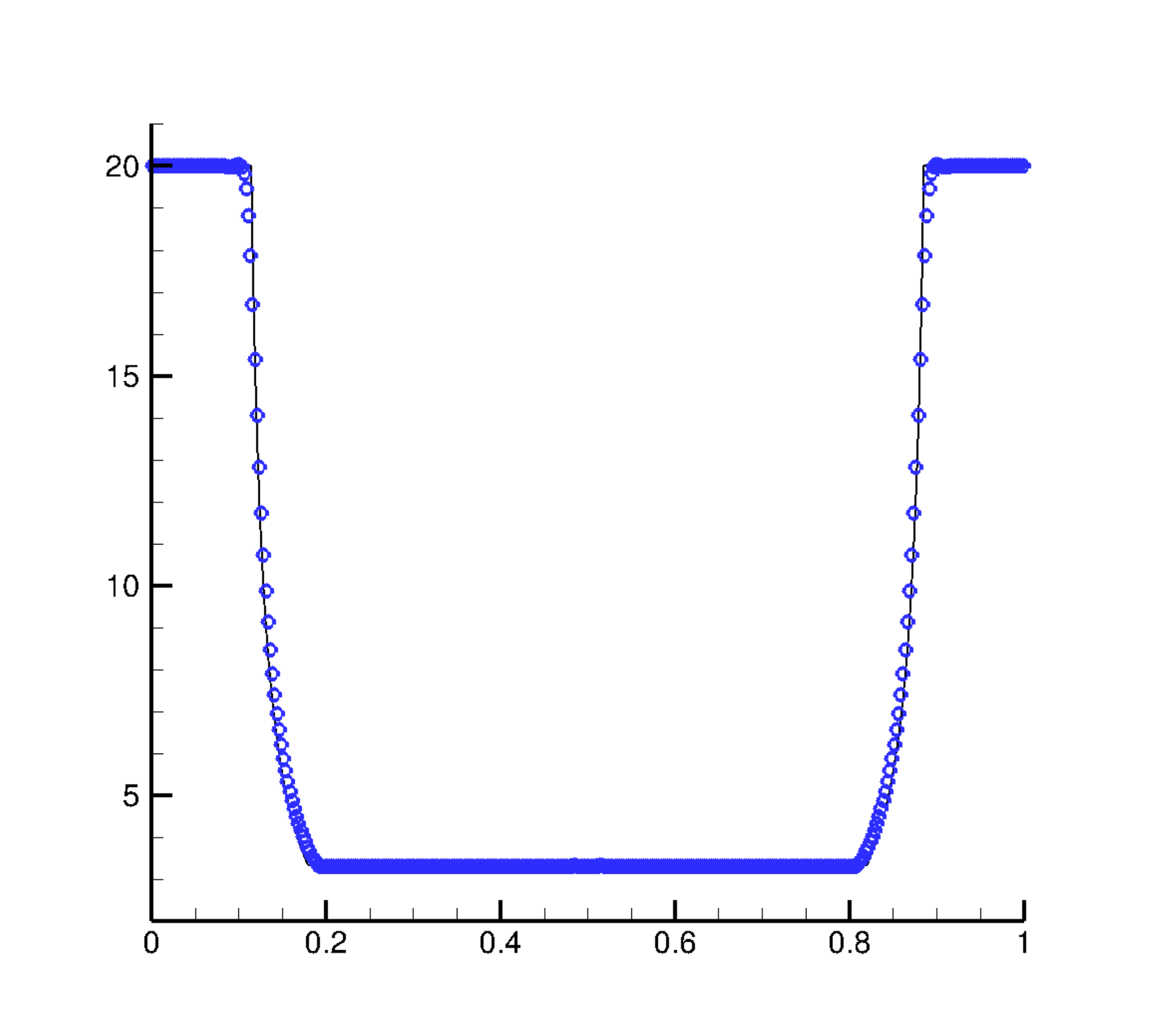}
  }
  \caption{Example \ref{ex:RP4}: Circles and
  the solid lines denote the numerical and  exact solutions, respectively. $N=400,~t=0.4$.}
  \label{fig:RP4}
\end{figure}

\begin{example}[Riemann problem 4]\label{ex:RP4}\rm
  The initial data of the fourth 1D Riemann problem are
  \begin{equation}
    (\rho,u,p)=\begin{cases}
      (1,~-0.7,~20),    &\quad x<0.5, \\
      (1,~0.7,~20), &\quad x>0.5.
    \end{cases}
  \end{equation}
\end{example}
The solution of this Riemann problem consists of a left-moving rarefaction wave,
a contact discontinuity, and a right-moving rarefaction wave.
The rest-mass density, the velocity, and the pressure at $t=0.4$ obtained by the
entropy stable scheme with $400$ are shown in Figure \ref{fig:RP4}.
It is seen that our entropy stable scheme can well capture the wave pattern, although
in the profile of the rest-mass density, there is undershoot similar to the results in \cite{yz1}.

\begin{figure}[!ht]
  \centering
  \subfigure[$\rho$]{
    \includegraphics[width=0.3\textwidth, trim=40 40 50 50, clip]{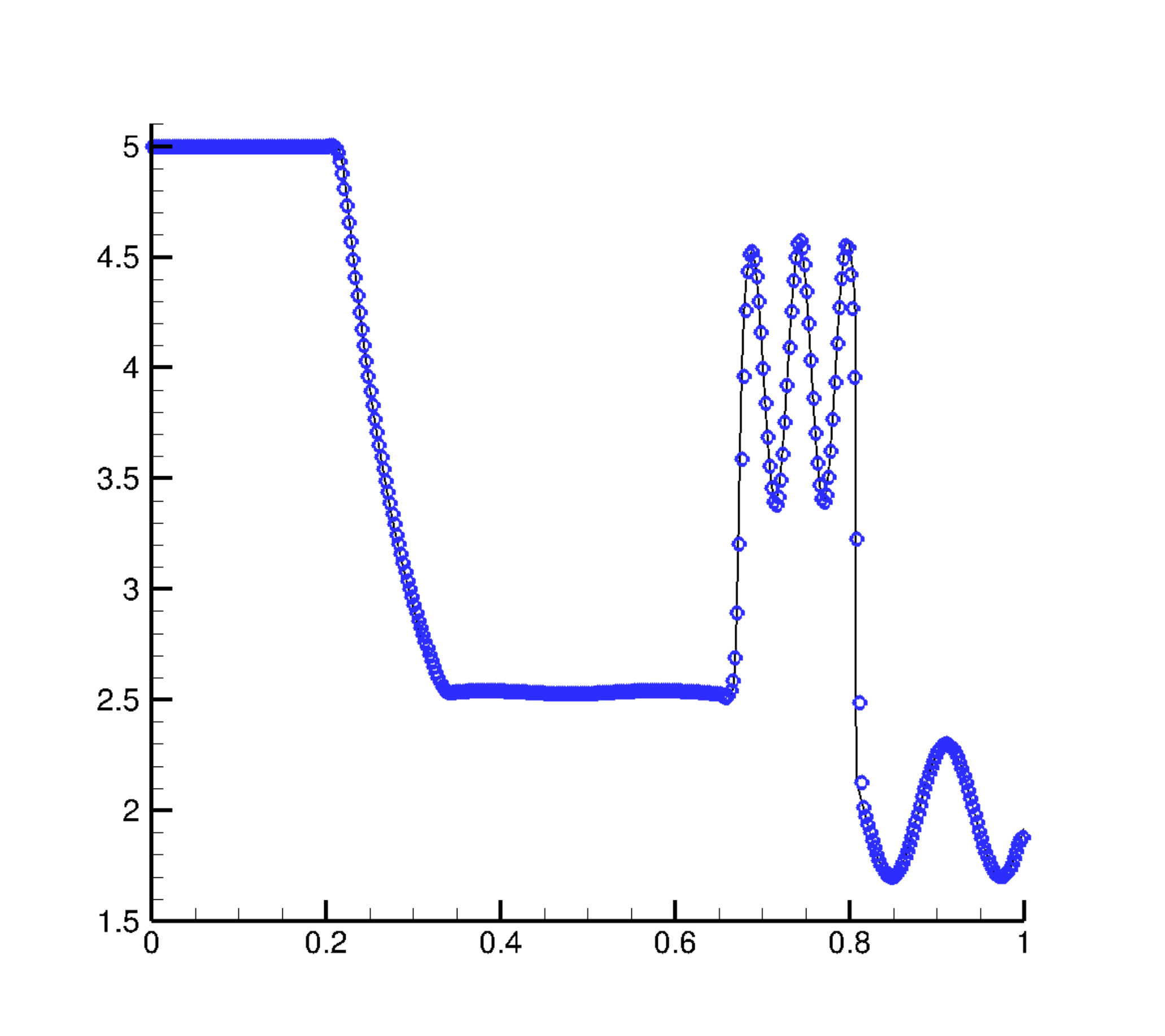}
  }
  \subfigure[$u$]{
    \includegraphics[width=0.3\textwidth, trim=40 40 50 50, clip]{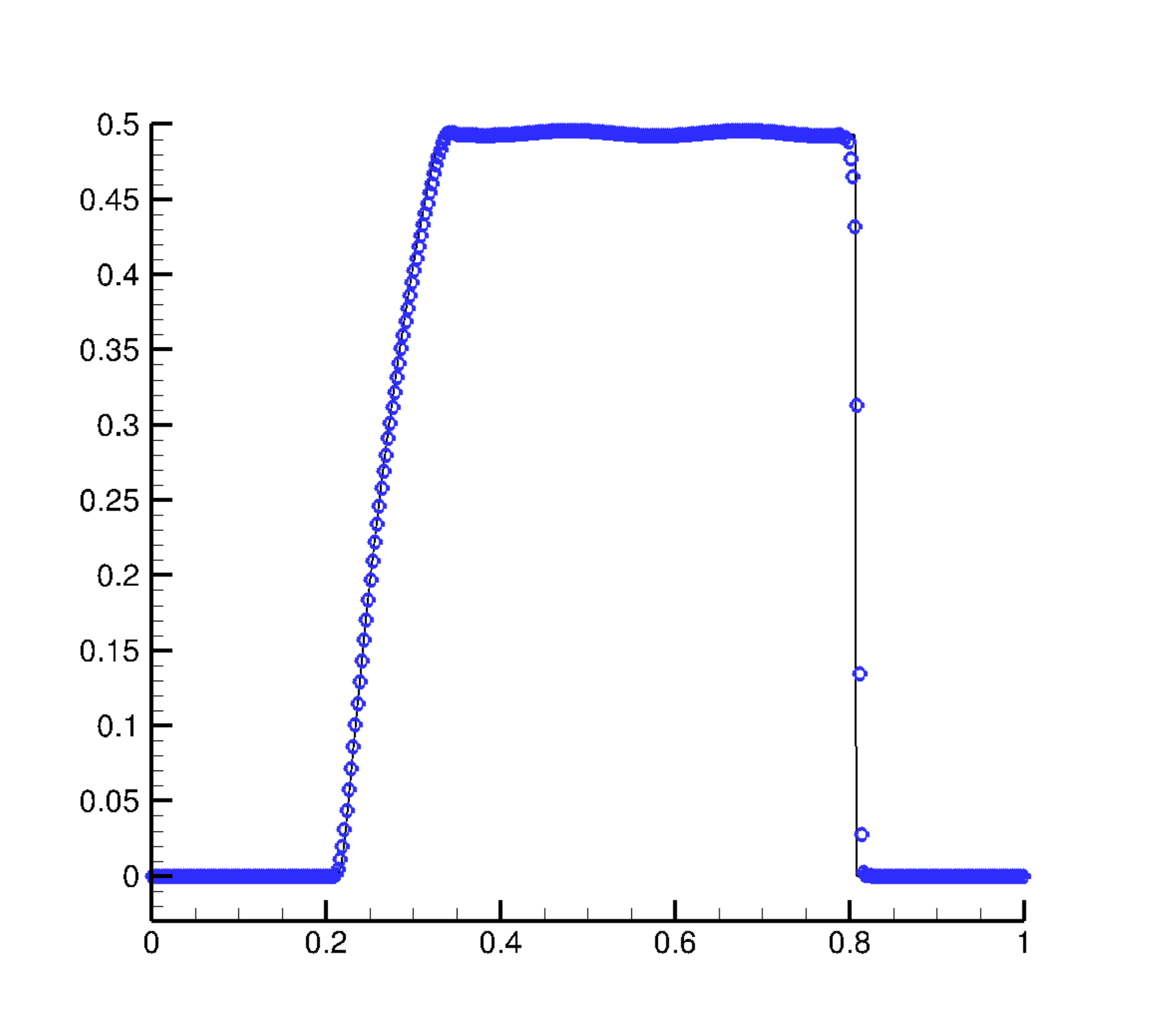}
  }
  \subfigure[$p$]{
    \includegraphics[width=0.3\textwidth, trim=40 40 50 50, clip]{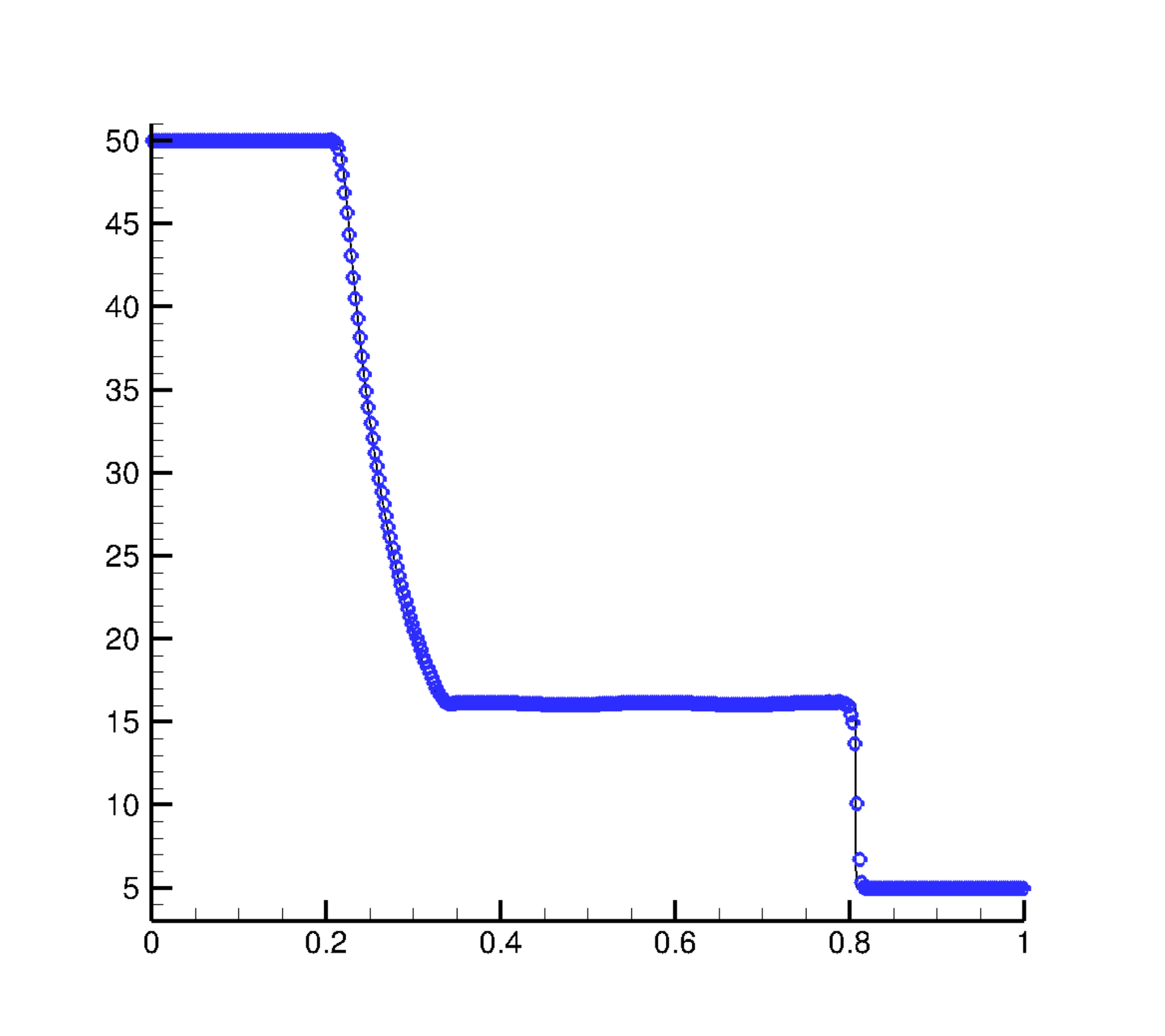}
  }
  \caption{Example \ref{ex:DP}: Circles and
  the solid lines denote the numerical and  reference solutions, respectively. $N=400,~t=0.35$.}
  \label{fig:DP}
\end{figure}

\begin{example}[Density perturbation problem]\label{ex:DP}\rm
	This is a more general Cauchy problem obtained
  by including a rest-mass density perturbation in the initial data of corresponding Riemann
  problem in order to test the ability of shock-capturing schemes
  to resolve small scale flow features, which may give a good indication
  of the numerical (artificial) viscosity of the scheme.
  The initial data are given by
  \begin{equation}
    (\rho,u,p)=\begin{cases}
      (5,~0,~50),    &\quad x<0.5, \\
      (2+0.3\sin(50x),~0,~5), &\quad x>0.5.
    \end{cases}
  \end{equation}
\end{example}
Figure \ref{fig:DP} plots the numerical results at $t = 0.35$ obtained by using our
entropy stable scheme with $400$. The reference solution is obtained by
using the first-order local Lax-Friedrichs scheme with $10000$ uniform cells.
We can see that the shock wave is moving into a sinusoidal rest-mass density
field, and  then some
smooth but complex structures are generated at the left when the shock wave
interacts with the sine wave; and our entropy stable scheme can capture them well.

\begin{figure}[!ht]
  \centering
  \subfigure[$\rho$]{
    \includegraphics[width=0.3\textwidth, trim=40 40 50 50, clip]{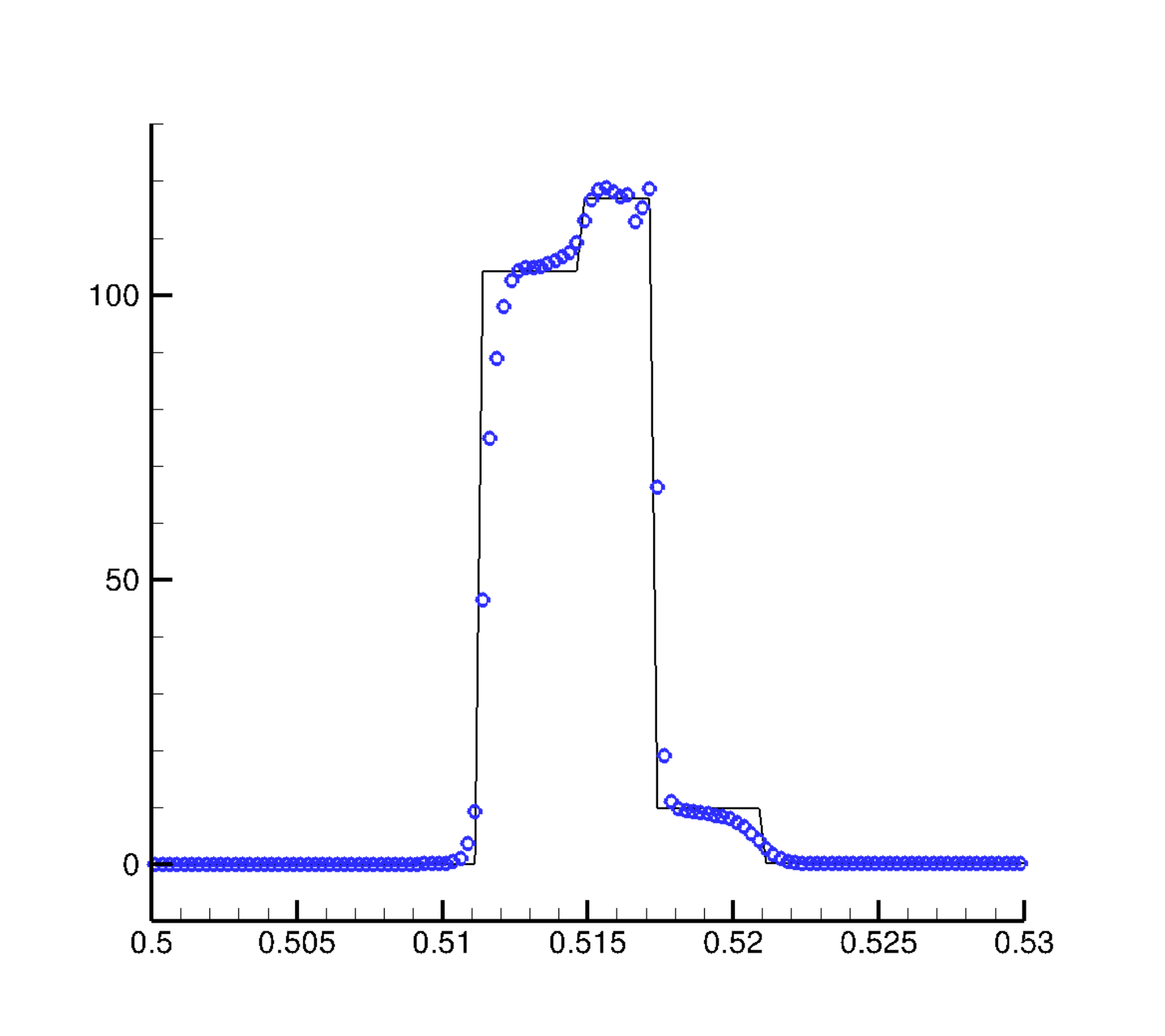}
  }
  \subfigure[$u$]{
    \includegraphics[width=0.3\textwidth, trim=40 40 50 50, clip]{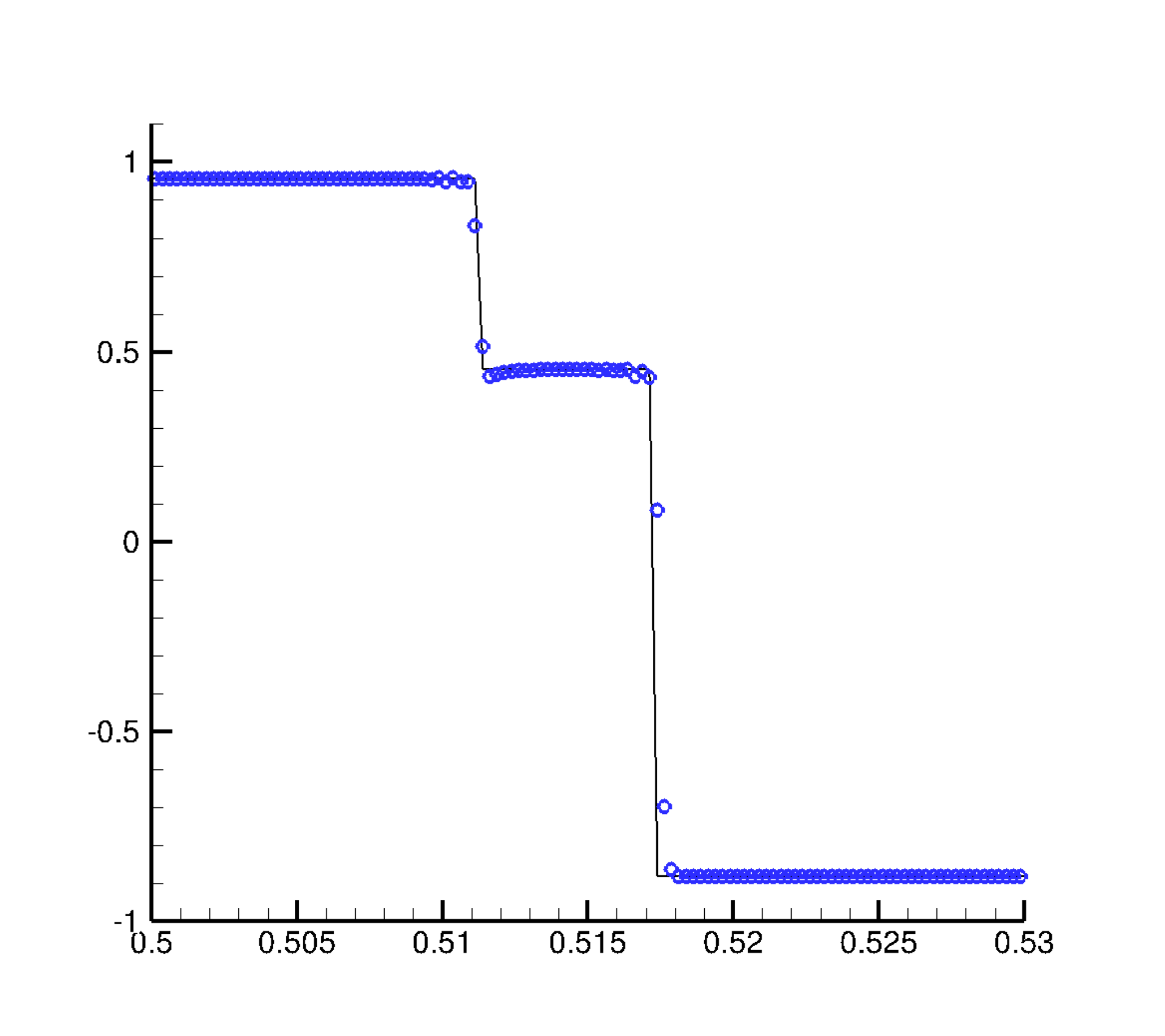}
  }
  \subfigure[$p$]{
    \includegraphics[width=0.3\textwidth, trim=40 40 50 50, clip]{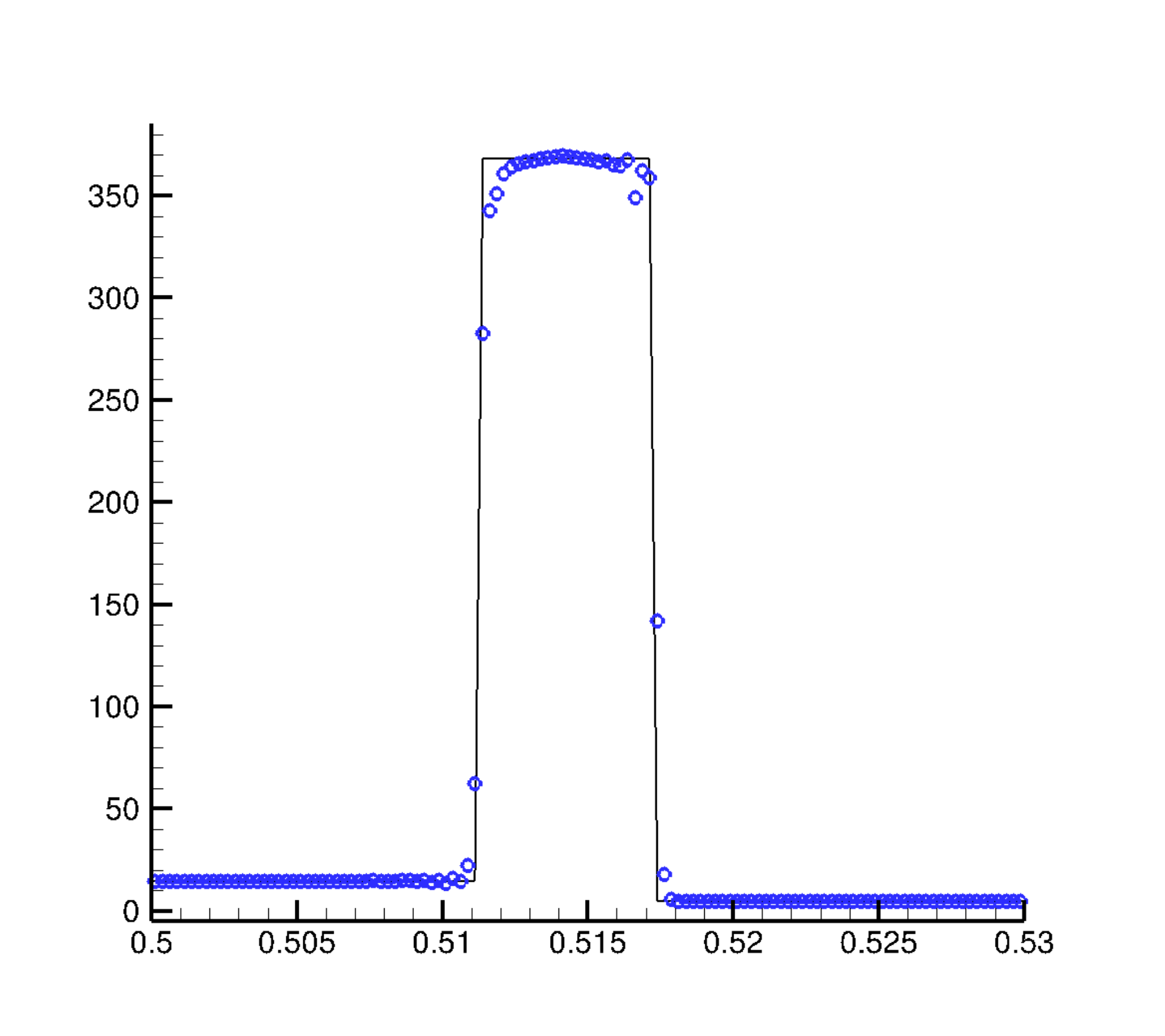}
  }
  \caption{Example \ref{ex:Blast}: Circles and
  the solid lines denote the numerical and  exact solutions, respectively. $N=4000,~t=0.43$.}
  \label{fig:Blast}
\end{figure}

\begin{example}[Blast wave interaction]\label{ex:Blast}\rm
  This test describes the collision of blast waves and is used to evaluate the
  performance of our entropy stable scheme for the flow with strong
  discontinuities. The initial condition is
  \begin{equation}
    (\rho,u,p)=\begin{cases}
      (1,~0,~10^3),    &\quad x<0.1, \\
      (1,~0,~10^{-2}),    &\quad 0.1<x<0.9, \\
      (1,~0,~10^2), &\quad x>0.9,
    \end{cases}
  \end{equation}
  with inflow and outflow boundary conditions and $\Gamma=1.4$.
\end{example}
Figure \ref{fig:Blast} gives close-up of the solutions at $t=0.43$ obtained by using
the entropy stable scheme with $4000$ uniform cells.
The solutions at this time within the interval $[0.5, 0.53]$
consists of two shock waves and two contact discontinuities.
It can be seen that our scheme can well resolve those
discontinuities and clearly capture the complex relativistic wave configuration,
except for slight overshoot and undershoot of the rest-mass density and the
pressure near $x=0.517$.

\subsection{Two-dimensional results}
\begin{example}[Accuracy test]\label{ex:acc2D}\rm
  It is a 2D relativistic isentropic vortex problem to test the accuracy,
  whose detailed construction can be found in \cite{Ling2019}.
  The initial rest-mass density and pressure are
  \begin{align*}
    &\rho(x,y)=(1-C_1 e^{1-r^2})^{\frac{1}{\Gamma-1}},\quad p=\rho^\Gamma,
  \end{align*}
  where
  \begin{align*}
    &C_1=\dfrac{(\Gamma-1)/\Gamma}{8\pi^2}\epsilon^2, \quad r=\sqrt{x_0^2+y_0^2},\\
    &x_0=x+\dfrac{\gamma-1}{2}(x+y), \quad y_0=y+\dfrac{\gamma-1}{2}(x+y),\quad \gamma = \dfrac{1}{\sqrt{1-w^2}},
  \end{align*}
  and the initial velocities are
  \begin{align*}
    &u=\dfrac{1}{1-{w(u_0+v_0)}/\sqrt{2}}\left[\dfrac{u_0}{\gamma}-\dfrac{w}{\sqrt{2}}+\dfrac{\gamma
    w^2}{2(\gamma+1)}(u_0+v_0)\right],\\
    &v=\dfrac{1}{1-{w(u_0+v_0)}/\sqrt{2}}\left[\dfrac{v_0}{\gamma}-\dfrac{w}{\sqrt{2}}+\dfrac{\gamma
    w^2}{2(\gamma+1)}(u_0+v_0)\right],
  \end{align*}
with
  \begin{align*}
    \quad (u_0,v_0)=(-y_0,x_0)f, \quad f=\sqrt{\dfrac{C_2}{1+C_2 r^2}},
    \quad C_2=\dfrac{2\Gamma C_1 e^{1-r^2}}{2\Gamma-1-\Gamma C_1 e^{1-r^2}}.
  \end{align*}
  This test describes a relativistic vortex moves with a constant speed of magnitude
  $w$ in $(-1,-1)$ direction.
\end{example}

\begin{figure}[!ht]
  \centering
  \subfigure[$\rho$]{
    \includegraphics[width=0.3\textwidth, trim=60 40 90 50, clip]{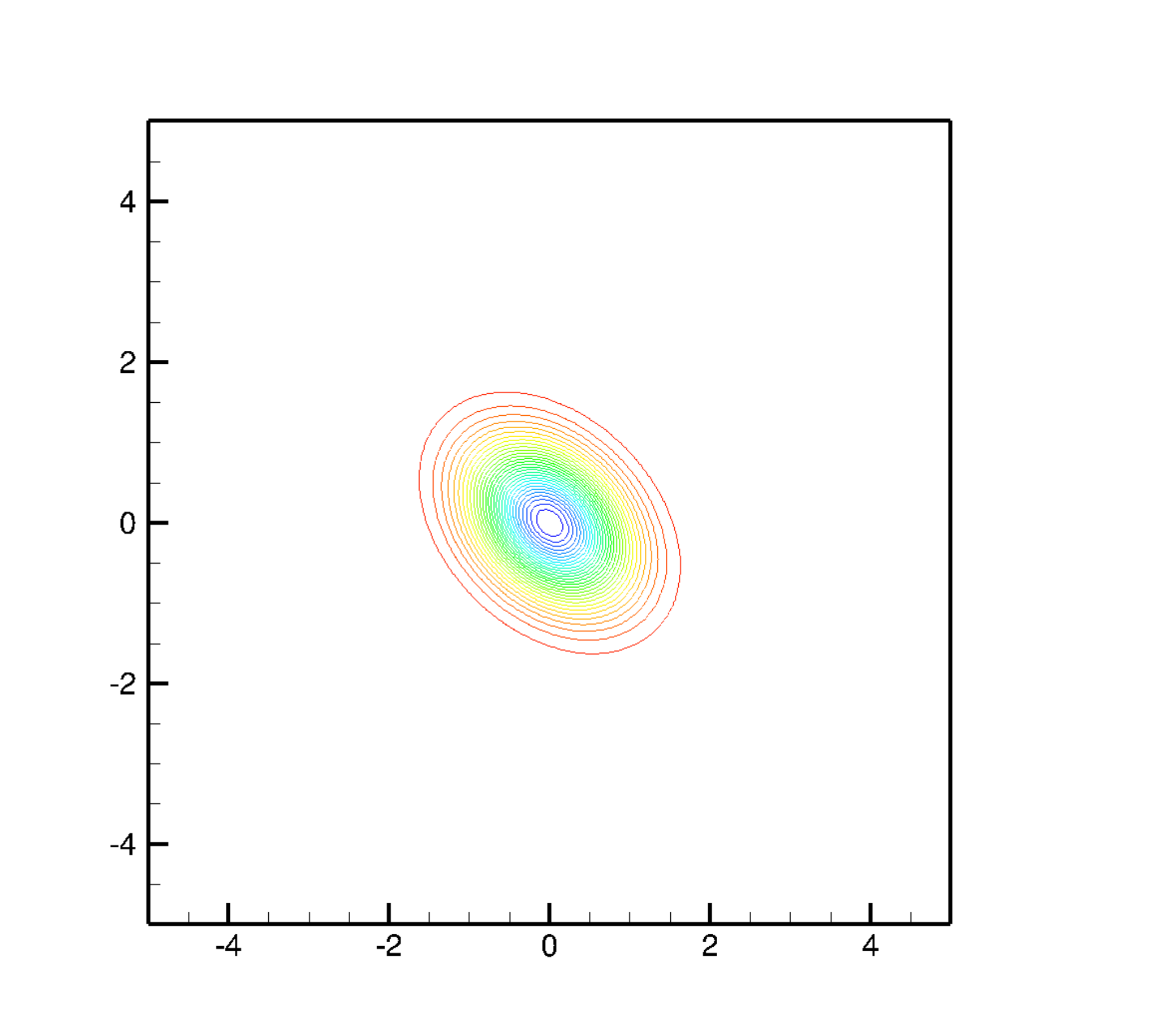}
  }
  \subfigure[$u$]{
    \includegraphics[width=0.3\textwidth, trim=60 40 90 50, clip]{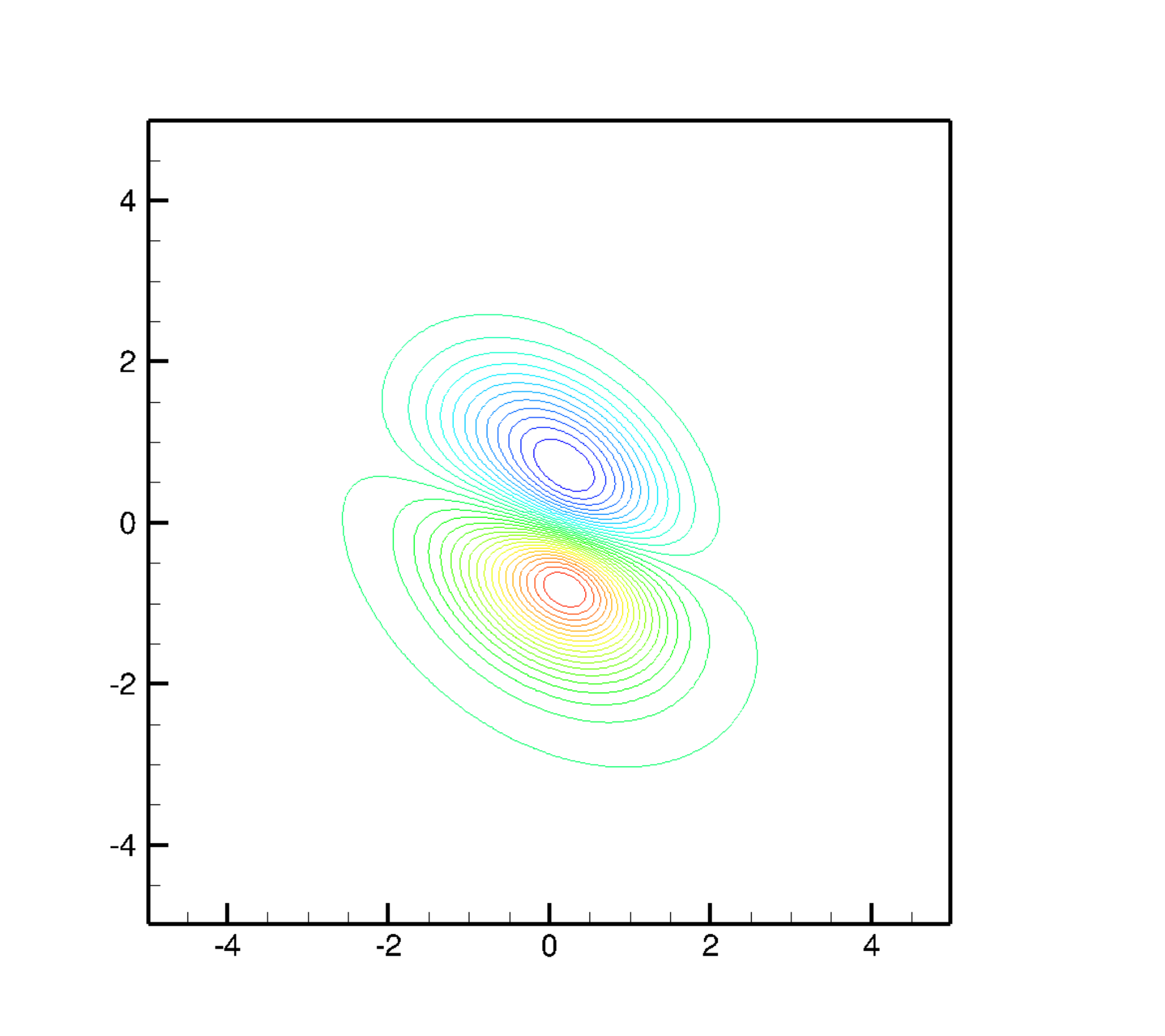}
  }
  \subfigure[$v$]{
    \includegraphics[width=0.3\textwidth, trim=60 40 90 50, clip]{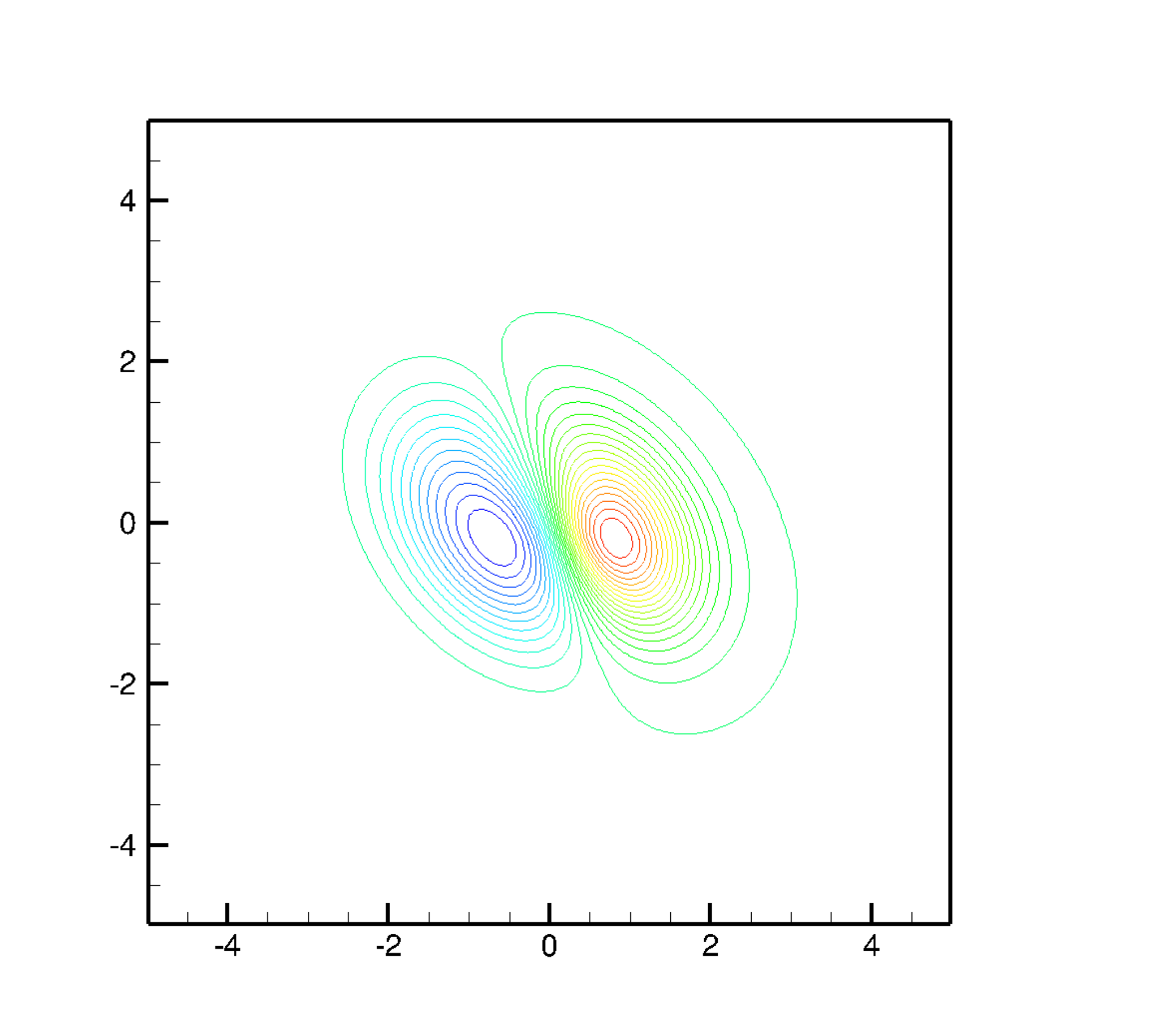}
  }
  \caption{Example \ref{ex:acc2D}: $30$ equally spaced contour lines.
  $N_x=N_y=320,~t=20$.}
  \label{fig:acc2D}
\end{figure}

\begin{figure}[!ht]
  \centering
  \includegraphics[width=0.50\textwidth, trim=20 40 50 50, clip]{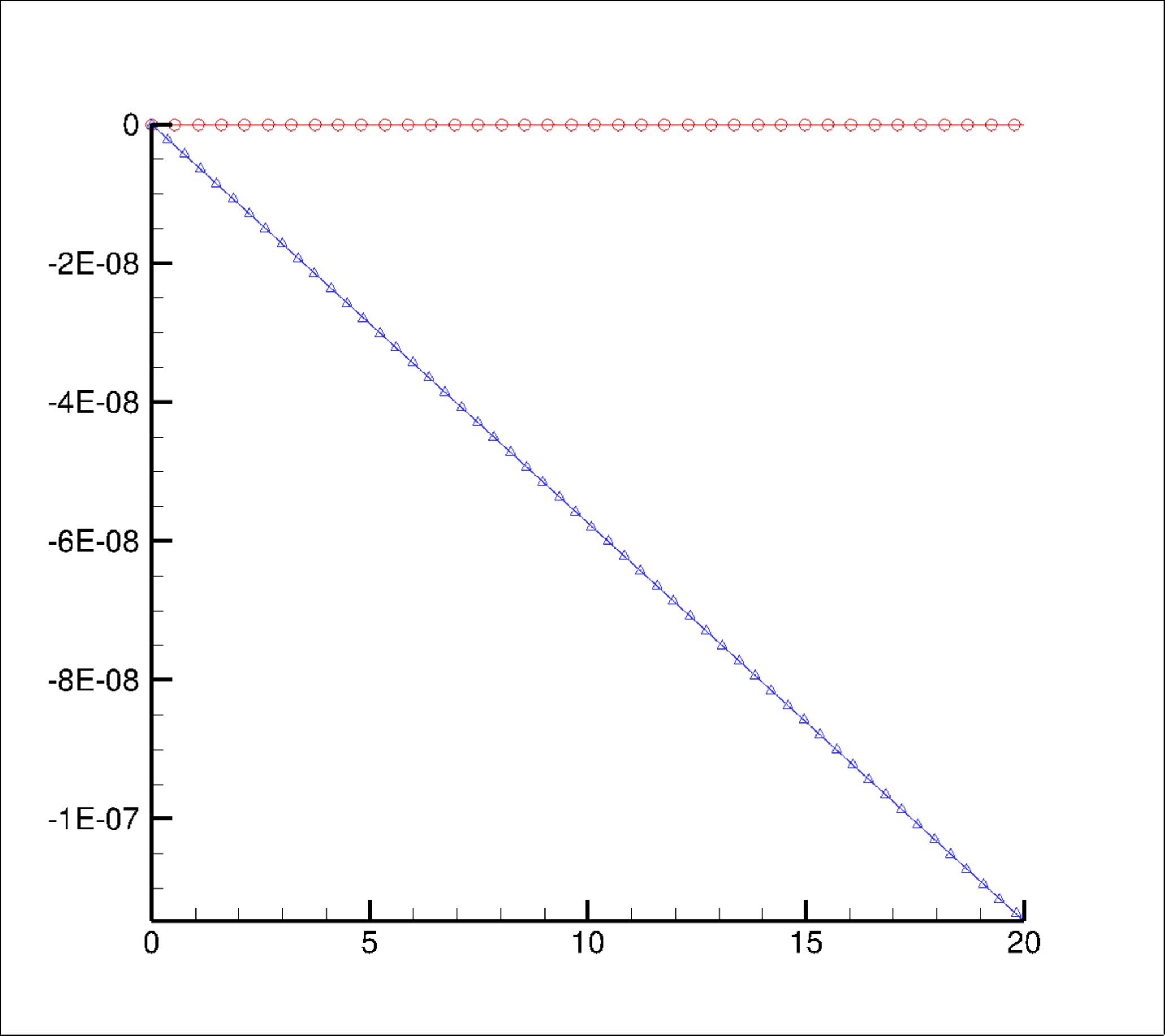}
  \caption{Example \ref{ex:acc2D}: The change of the total entropy with respect to $t$.
    Circles and deltas denote the results obtained by the entropy conservative scheme and
  the entropy stable scheme, respectively.   $N_x=N_y=320$.}
  \label{fig:decay}
\end{figure}

The test is computed in the domain $[-5,5]^2$ with $w=0.5\sqrt{2},\epsilon=5$, and periodic boundary
conditions. The output time is $t=20$ so that the vortex travels and returns to the original position after a period.
Table \ref{tab:acc2D} lists the errors of the rest-mass density and orders of convergence.
It can be clearly seen that our entropy stable scheme achieves fifth-order accuracy.
Figure \ref{fig:acc2D} plots the contours of the rest-mass density, and the
velocities with $30$ equally spaced contour lines.
The results show that due to the Lorentz contraction, the vortex becomes
elliptic and the velocity in $x$- (resp. $y$-) direction is not symmetric respect to
$y=0$ (resp. $x=0$) anymore.
Figure \ref{fig:decay} presents
the change of the total entropy
$\sum_{i,j}\eta(\bU_{i,j})\Delta x\Delta y$
with respect to the time obtained by the entropy conservative scheme
and the entropy stable scheme with $N_x=N_y=320$.
We can see that for the entropy conservative scheme,
the total entropy almost keeps conservative
and for the entropy stable scheme, the total entropy decays as expected.

\begin{table}[!ht]
  \centering
  \begin{tabular}{r|cc|cc|cc} \hline
 $N$& $\ell^1$ error & order & $\ell^2$ error & order & $\ell^\infty$ error & order \\ \hline
 20 & 1.704e-02 &  -   & 4.982e-02 &  -   & 4.276e-01 &  -   \\
 40 & 2.886e-03 & 2.56 & 8.947e-03 & 2.48 & 7.352e-02 & 2.54 \\
 80 & 1.781e-04 & 4.02 & 6.750e-04 & 3.73 & 1.300e-02 & 2.50 \\
160 & 4.973e-06 & 5.16 & 2.962e-05 & 4.51 & 9.425e-04 & 3.79 \\
320 & 1.026e-07 & 5.60 & 8.240e-07 & 5.17 & 3.048e-05 & 4.95 \\
\hline
  \end{tabular}
  \caption{Example \ref{ex:acc2D}: Errors and orders of convergence in the rest-mass density
  at $t = 20$.}
  \label{tab:acc2D}
\end{table}

\begin{example}[Riemann problem 1]\label{ex:2DRP1}\rm
  The initial data are
  \begin{align*}
    (\rho,u,v,p)=\begin{cases}
      (0.5,~0.5,-0.5,~5), &\quad x>0.5,~y>0.5,\\
      (1,~0.5,~0.5,~5),    &\quad x<0.5,~y>0.5,\\
      (3,-0.5,~0.5,~5),   &\quad x<0.5,~y<0.5,\\
      (1.5,-0.5,-0.5,~5),&\quad x>0.5,~y<0.5.
    \end{cases}
  \end{align*}
  It describes the interaction of four contact discontinuities (vortex sheets)
  with the same sign (the negative sign).
\end{example}

\begin{figure}[!ht]
  \centering
  \includegraphics[width=0.45\textwidth, trim=40 40 90 50, clip]{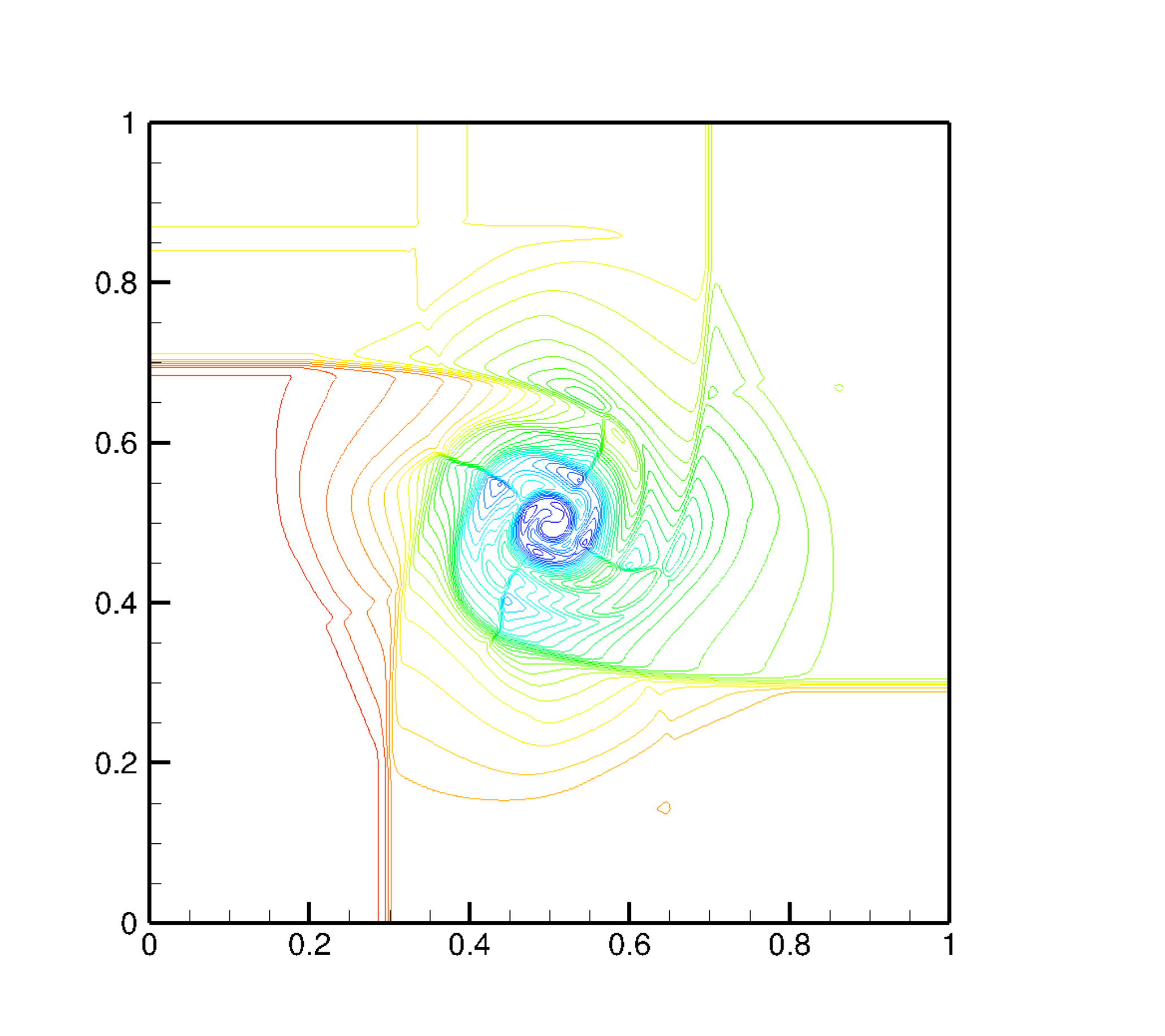}
  \includegraphics[width=0.45\textwidth, trim=40 40 90 50, clip]{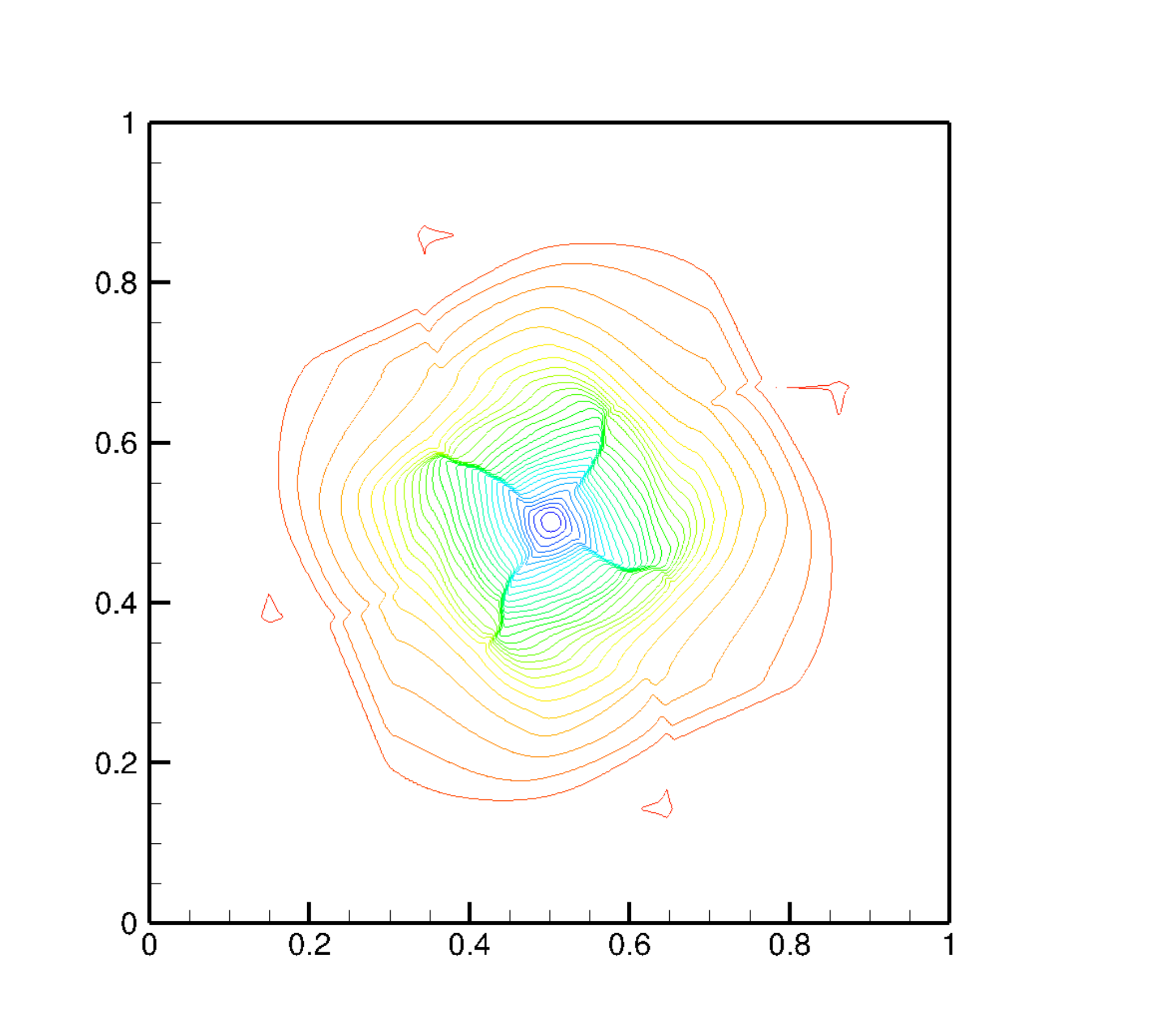}
  \caption{Example \ref{ex:2DRP1}: $N_x=N_y=400$, left: $\ln \rho$, right: $\ln p$, $30$ equally spaced contour lines.}
  \label{fig:2DRP1}
\end{figure}

Figure \ref{fig:2DRP1} shows the contours of the rest-mass density and pressure
logarithms with $30$ equally spaced contour lines.
We can see that the four initial vortex sheets interact each other to
form a spiral with the low rest-mass density around the center of the domain as time
increases, which is the typical cavitation phenomenon in gas dynamics.

\begin{figure}[!ht]
  \centering
  \includegraphics[width=0.45\textwidth, trim=40 40 90 50, clip]{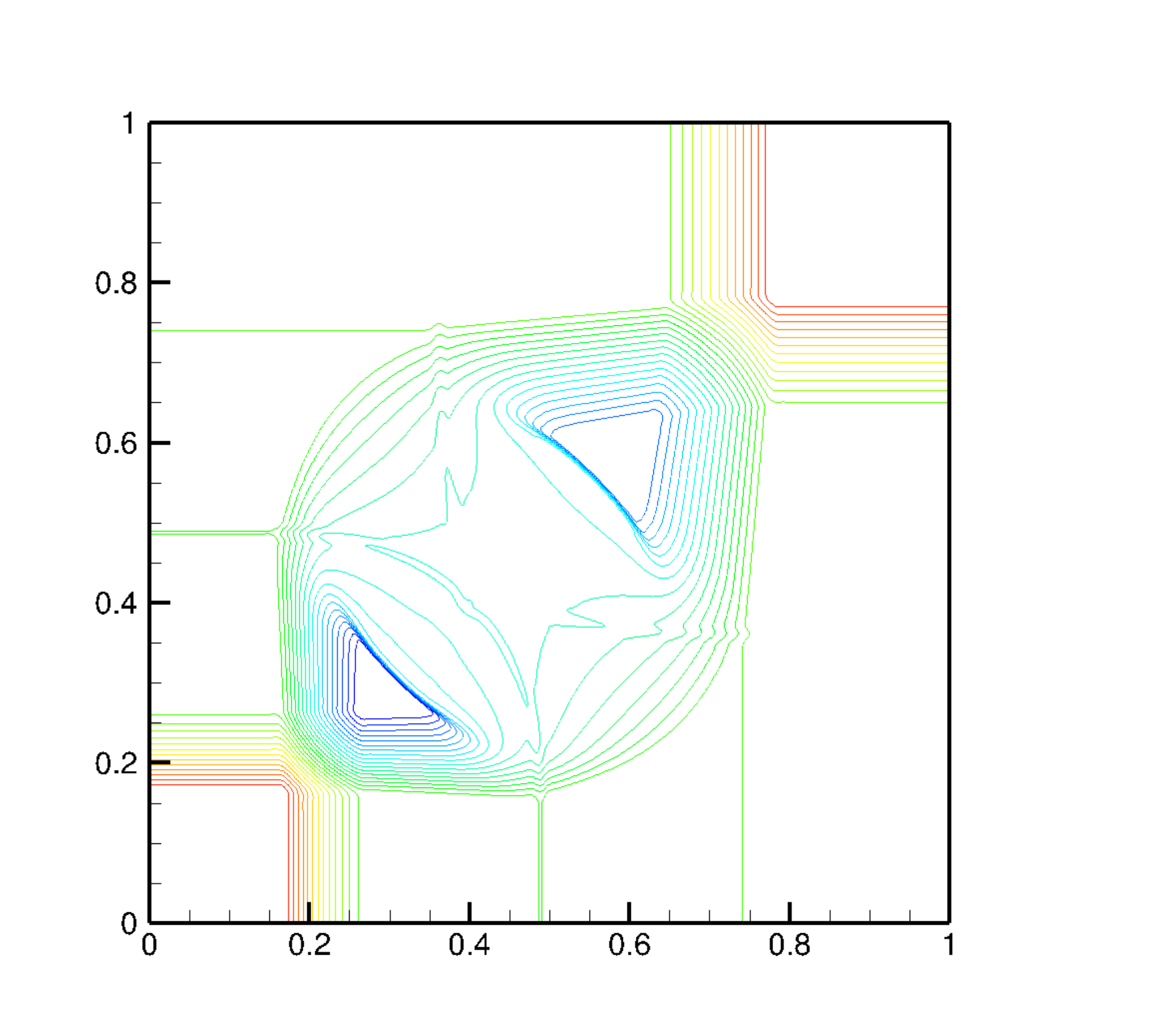}
  \includegraphics[width=0.45\textwidth, trim=40 40 90 50, clip]{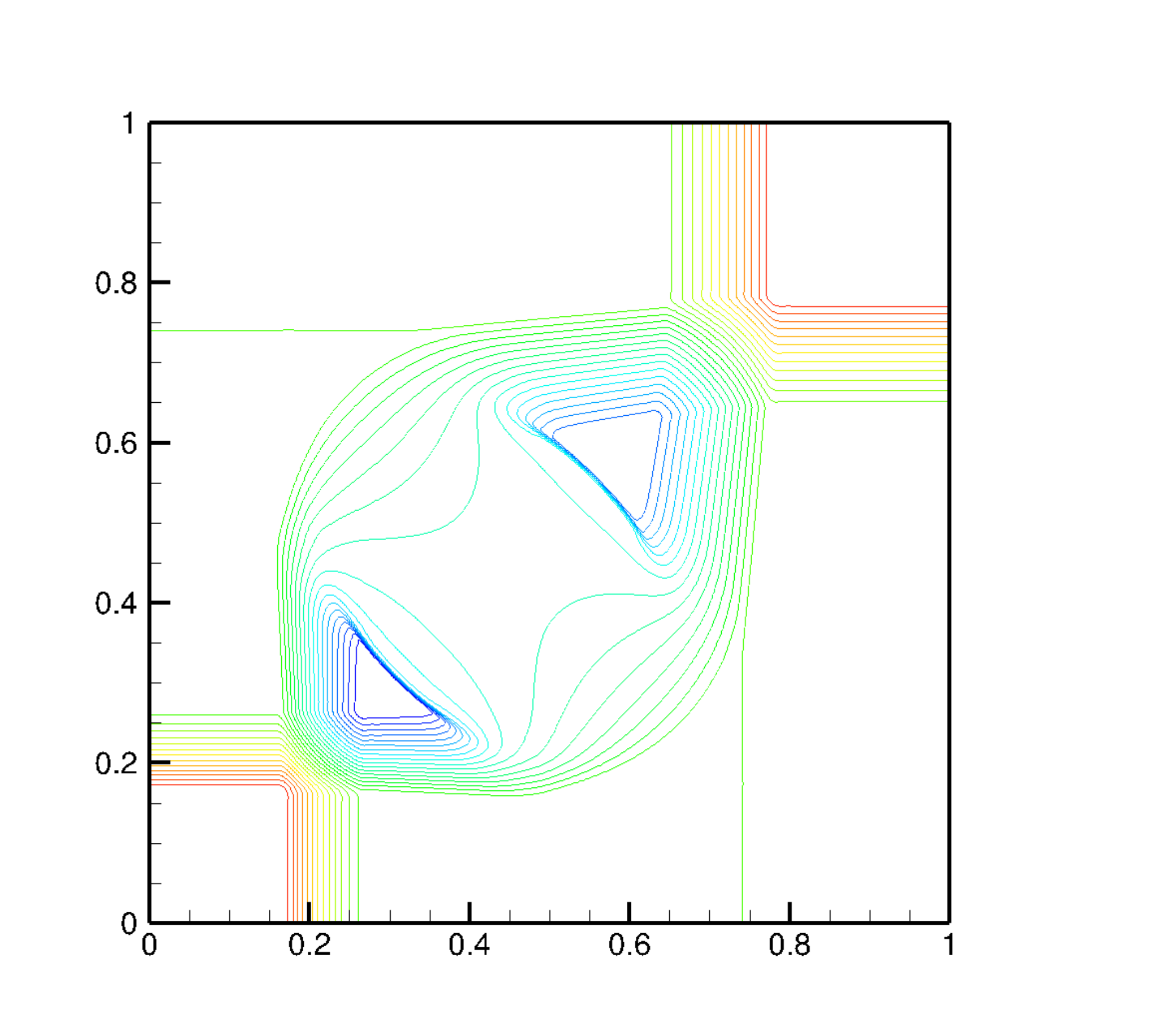}
  \caption{Example \ref{ex:2DRP2}: $N_x=N_y=400$, left: $\ln \rho$, right: $\ln p$, $30$ equally spaced contour lines.}
  \label{fig:2DRP2}
\end{figure}

\begin{example}[Riemann problem 2]\label{ex:2DRP2}\rm
  The initial data are
  \begin{align*}
    (\rho,u,v,p)=\begin{cases}
      (1,~0,~0,~1), &\quad x>0.5,~y>0.5,\\
      (0.5771,-0.3529,~0,~0.4),    &\quad x<0.5,~y>0.5,\\
      (1,-0.3529,-0.3529,~1),   &\quad x<0.5,~y<0.5,\\
      (0.5771,~0,-0.3529,~0.4),&\quad x>0.5,~y<0.5,
    \end{cases}
  \end{align*}
  which is about the interaction of four rarefaction waves.
\end{example}
Figure \ref{fig:2DRP2} plots the contours of the rest-mass density and pressure logarithms with $30$ equally spaced contour lines.
The results show that those four
initial discontinuities first evolve as four rarefaction waves and then interact each other and form two
(almost parallel) curved shock waves perpendicular to the line $y=x$ as time increases.

\begin{figure}[!ht]
  \centering
  \includegraphics[width=0.45\textwidth, trim=40 40 90 50, clip]{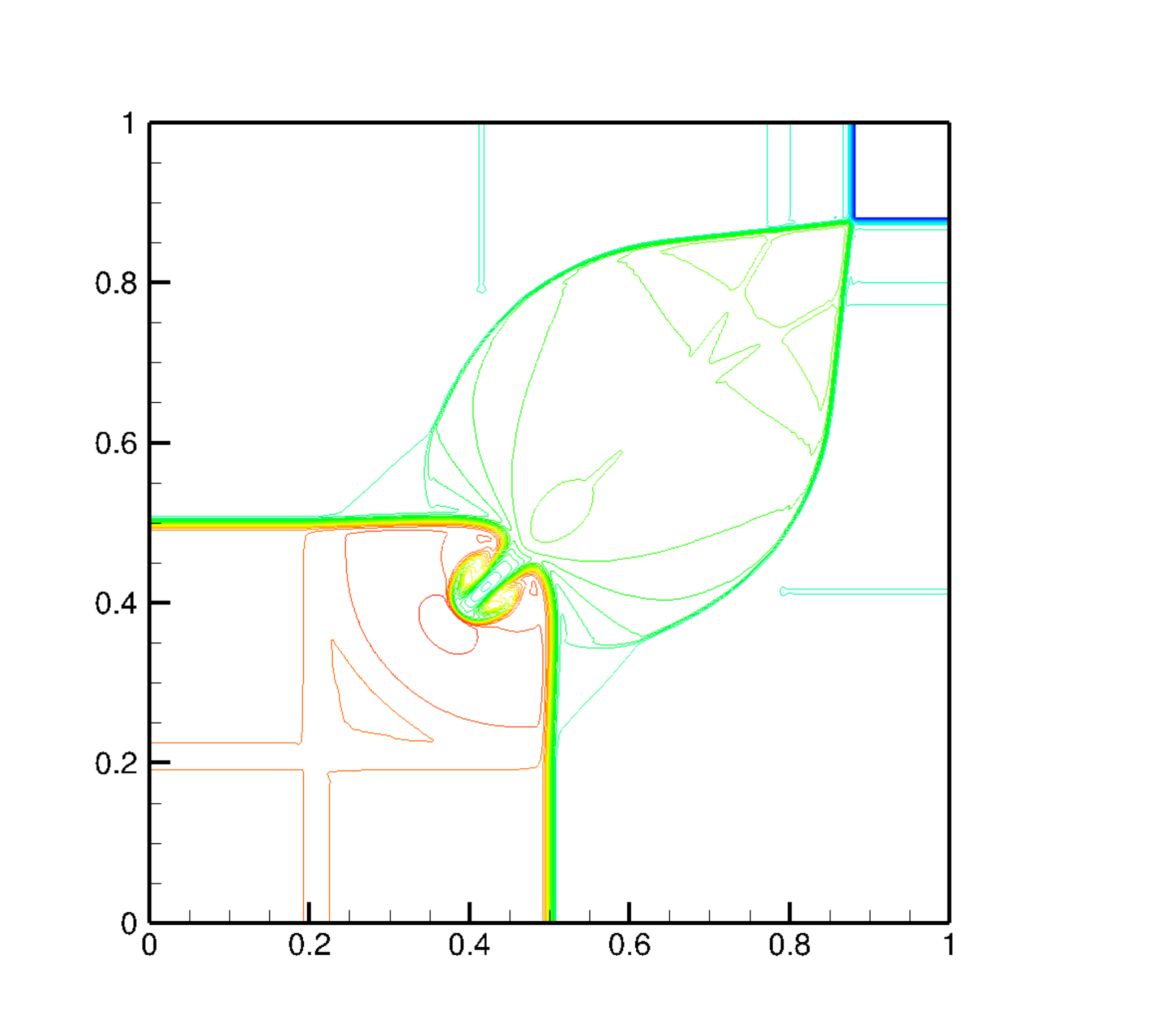}
  \includegraphics[width=0.45\textwidth, trim=40 40 90 50, clip]{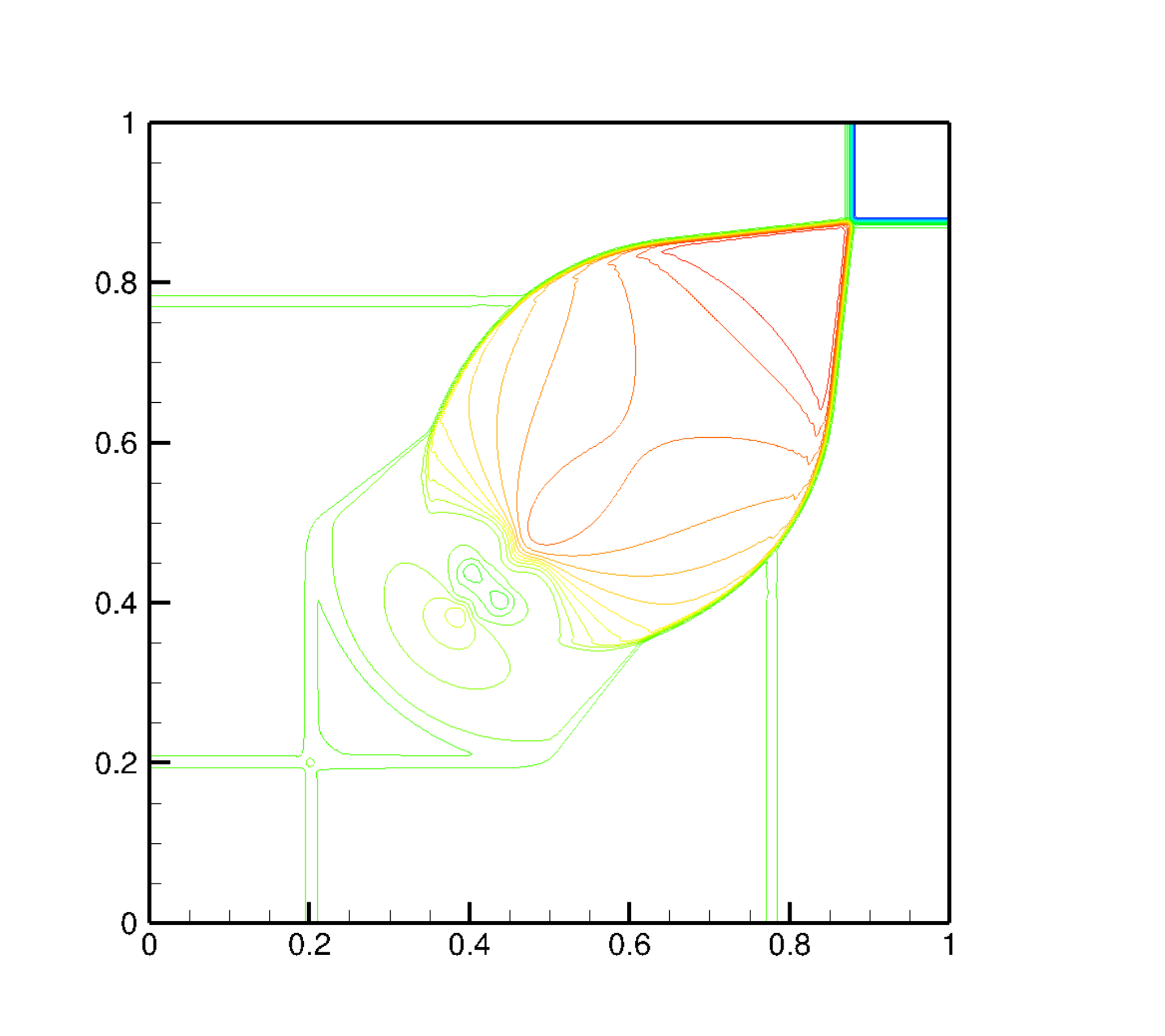}
  \caption{Example \ref{ex:2DRP3}: $N_x=N_y=400$, left: $\ln \rho$, right: $\ln p$, $30$ equally spaced contour lines.}
  \label{fig:2DRP3}
\end{figure}

\begin{example}[Riemann problem 3]\label{ex:2DRP3}\rm
  The initial data are
  \begin{align*}
    (\rho,u,v,p)=\begin{cases}
      (0.035145216124503,~0,~0,~0.162931056509027), &\quad x>0.5,~y>0.5,\\
      (0.1,~0.7,~0,~1),    &\quad x<0.5,~y>0.5,\\
      (0.5,~0,~0,~1),   &\quad x<0.5,~y<0.5,\\
      (0.1,~0,~0.7,~1),&\quad x>0.5,~y<0.5,
    \end{cases}
  \end{align*}
  where the left and bottom discontinuities are two contact discontinuities and
  the top and right are two shock waves.
\end{example}
We show the contours of the rest-mass density and pressure logarithms with $30$
equally spaced contour lines in Figure \ref{fig:2DRP3}.
The initial discontinuities interact each other,
and form a ``mushroom cloud'' around the point $(0.5, 0.5)$.

\begin{figure}[!ht]
  \centering
  \includegraphics[width=0.80\textwidth, trim=50 40 50 350, clip]{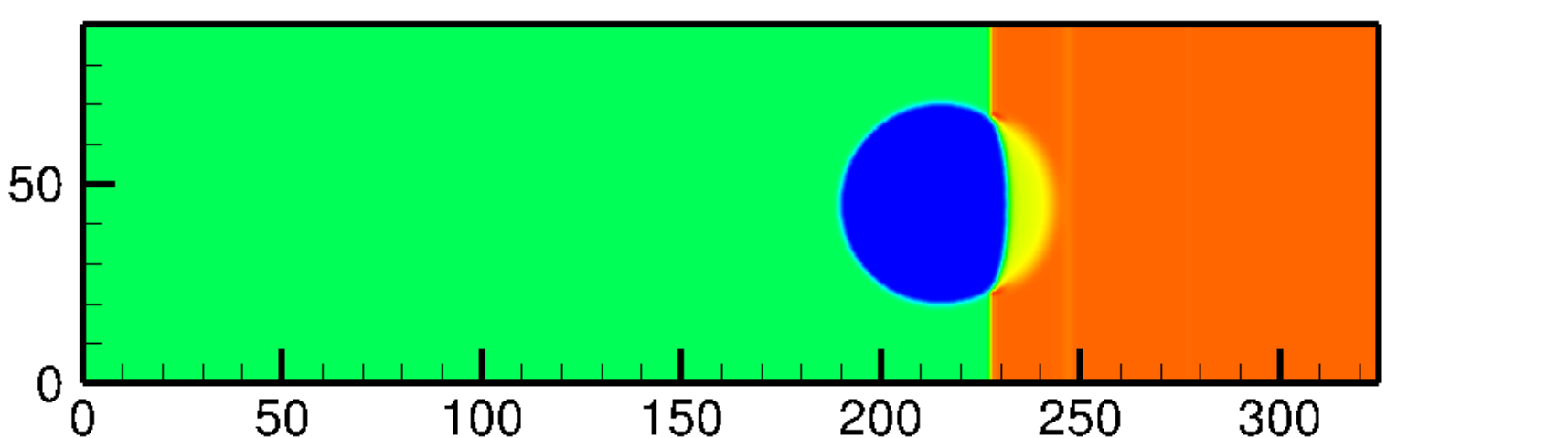}
  \includegraphics[width=0.80\textwidth, trim=50 40 50 350, clip]{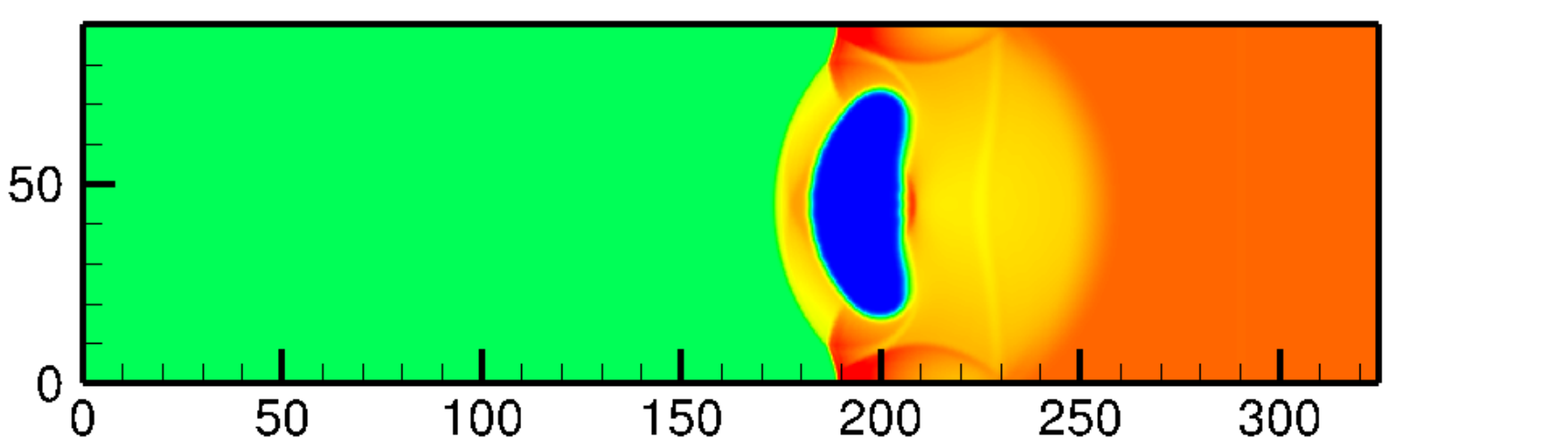}
  \includegraphics[width=0.80\textwidth, trim=50 40 50 350, clip]{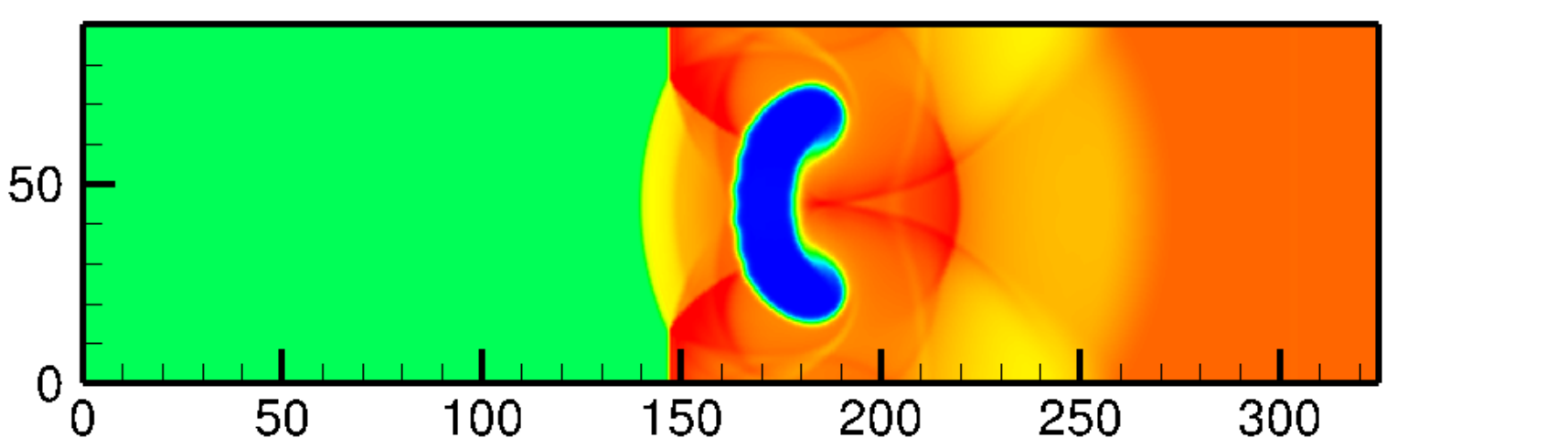}
  \includegraphics[width=0.80\textwidth, trim=50 40 50 350, clip]{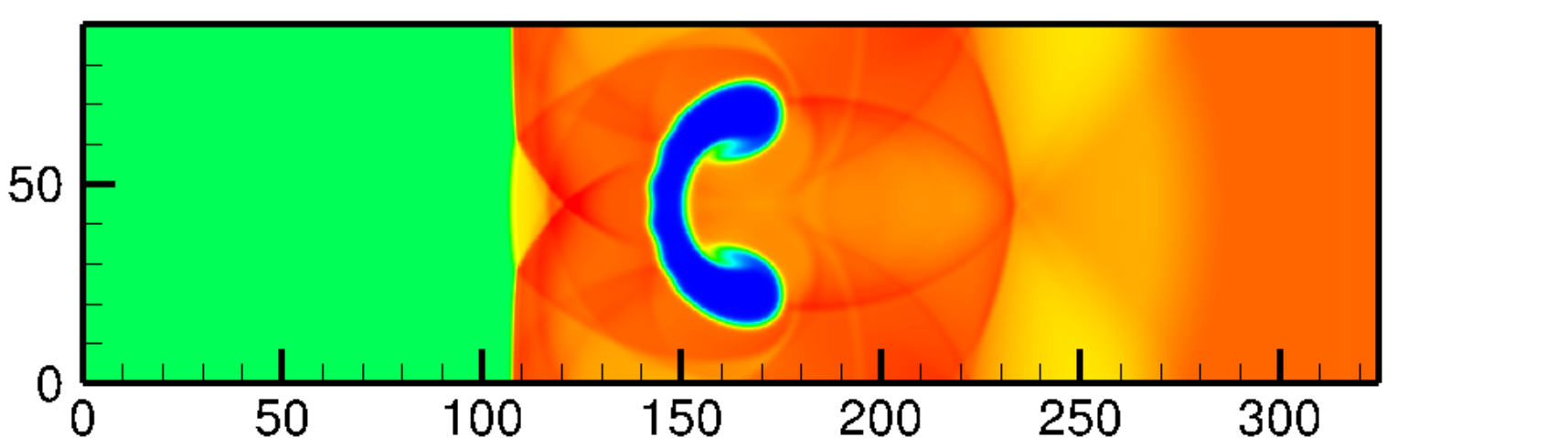}
  \includegraphics[width=0.80\textwidth, trim=50 40 50 350, clip]{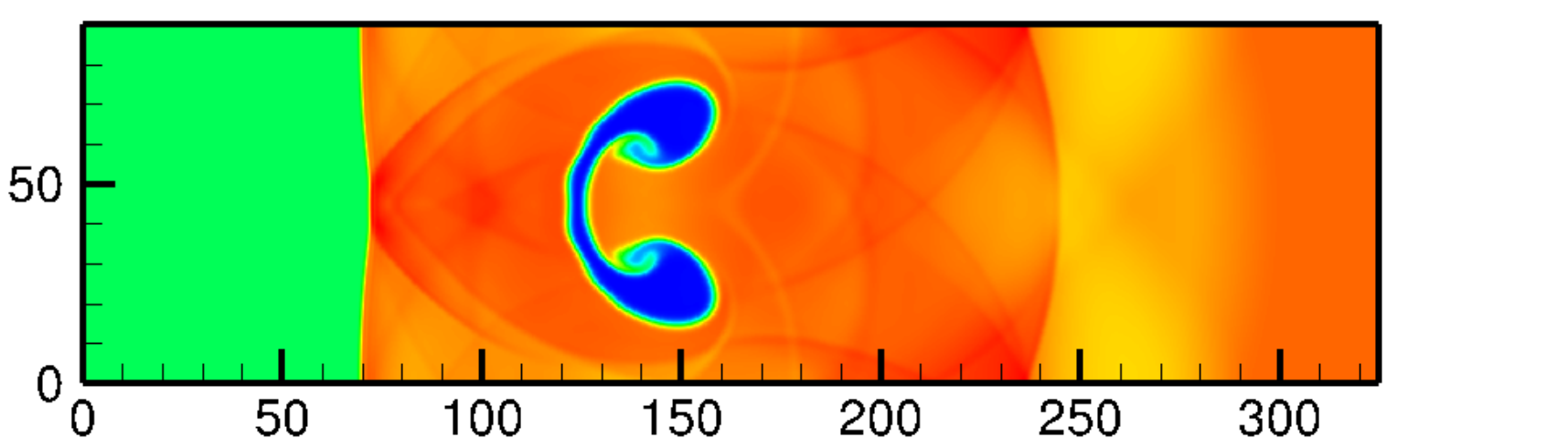}
  \caption{The first problem of Example \ref{ex:SB}: the schlieren images of $\rho$ at  $t=90,180,270,360,450$ from top to
    bottom. }
  \label{fig:SB1sch}
\end{figure}

\begin{figure}[!ht]
  \centering
  \includegraphics[width=0.80\textwidth, trim=50 40 50 350, clip]{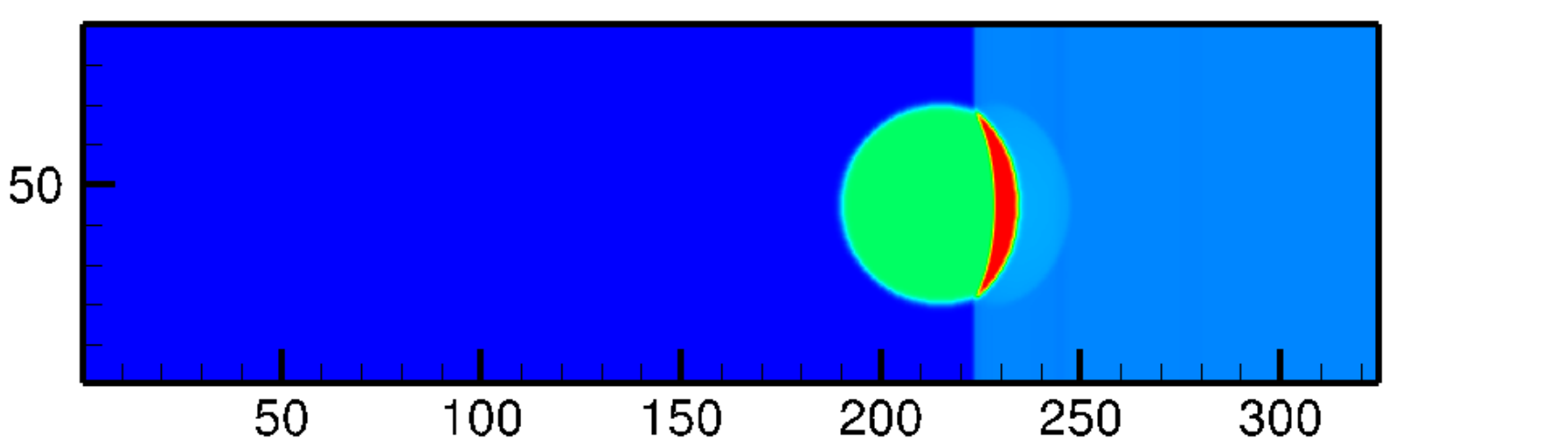}
  \includegraphics[width=0.80\textwidth, trim=50 40 50 350, clip]{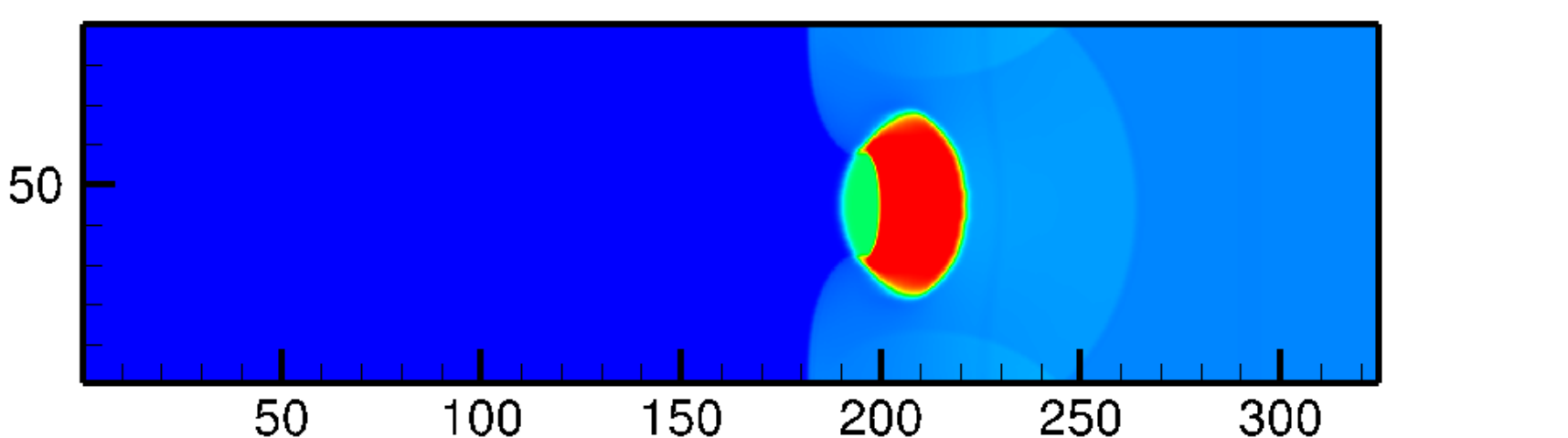}
  \includegraphics[width=0.80\textwidth, trim=50 40 50 350, clip]{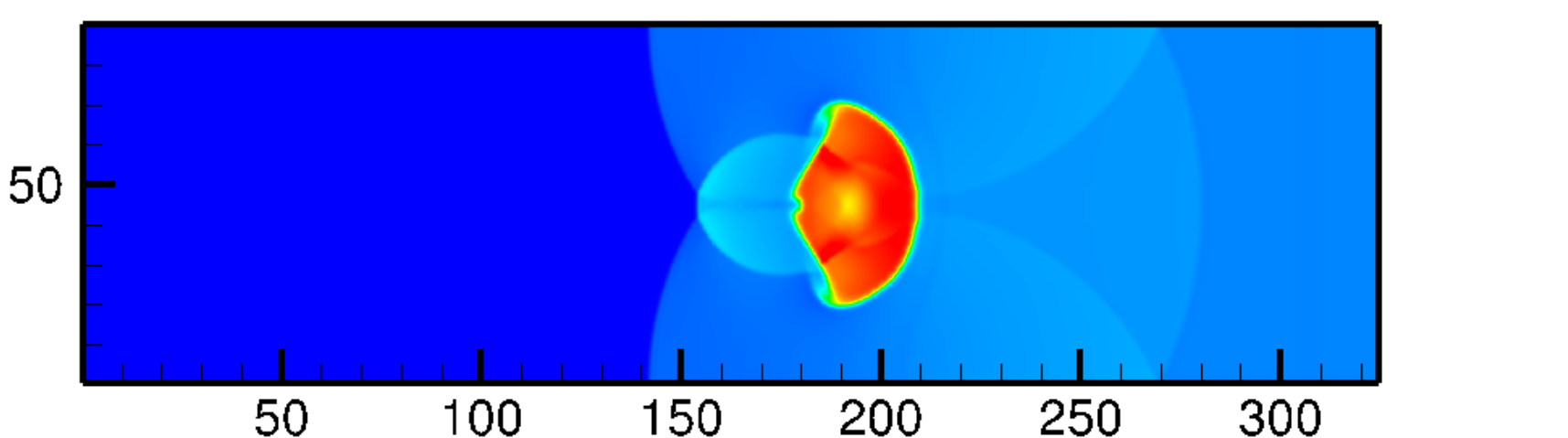}
  \includegraphics[width=0.80\textwidth, trim=50 40 50 350, clip]{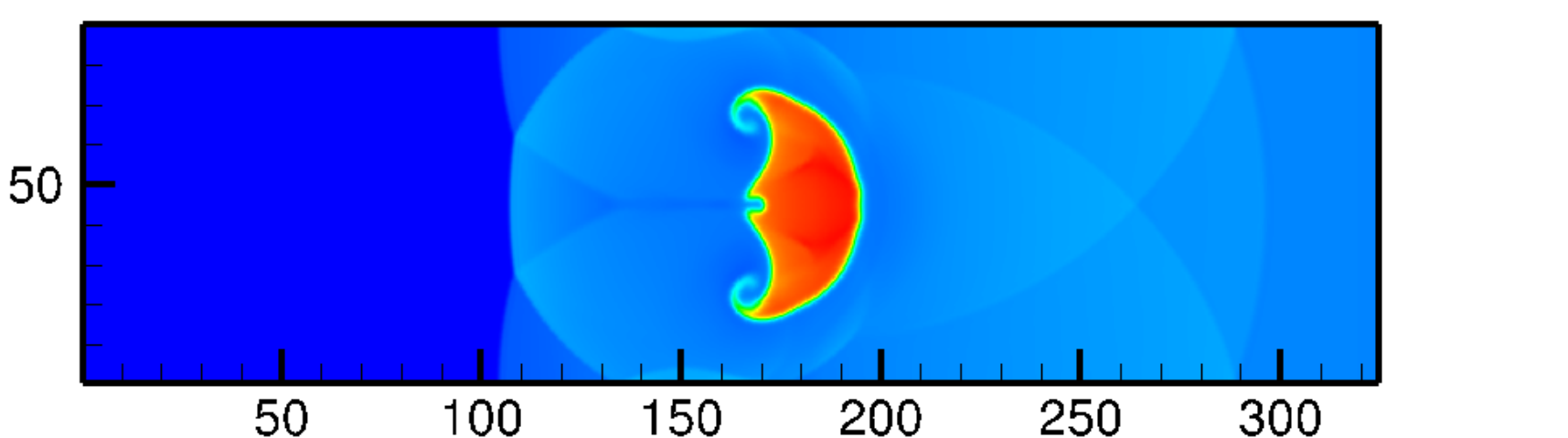}
  \includegraphics[width=0.80\textwidth, trim=50 40 50 350, clip]{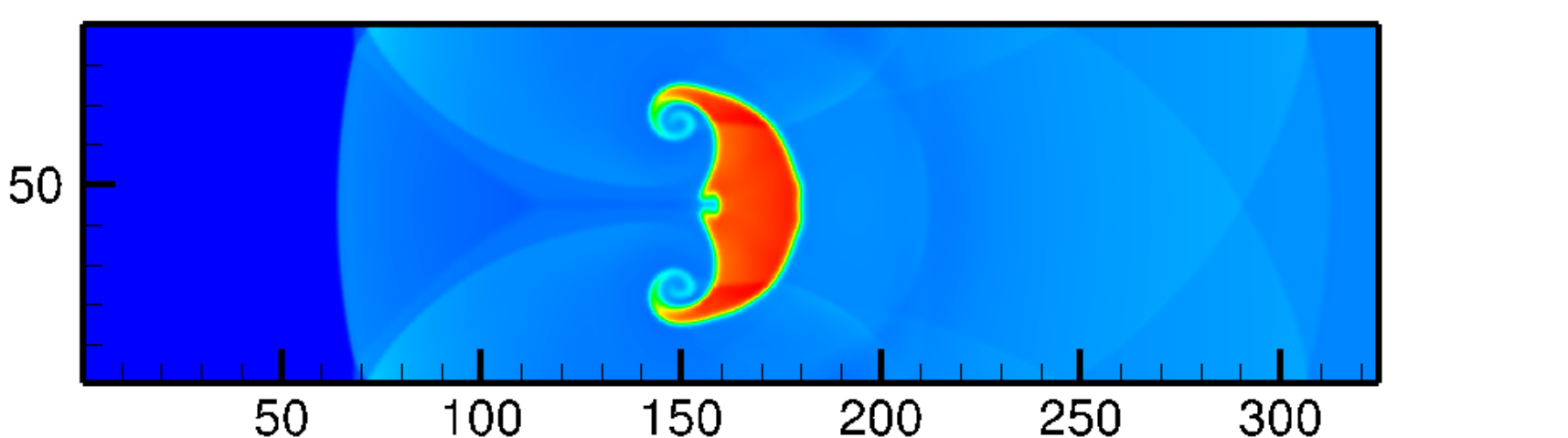}
  \caption{The second problem of Example \ref{ex:SB}: the schlieren images of $\rho$ at $t=100,200,300,400,500$ from top to bottom. }
  \label{fig:SB2sch}
\end{figure}

\begin{example}[Shock-bubble interaction problems]\label{ex:SB}\rm
  This example considers two shock-bubble interaction problems within the computational domain
  $[0,325]\times[-45,45]$. The detailed setup can be found in \cite{he1}.
  For the first problem, the initial left and right states of the planar shock
  wave moving left are
  \begin{equation*}
    (\rho,u,v,p)=\begin{cases}
      (1,~0,~0,~0.05),&\quad x<265, \\
      (1.865225080631180,-0.196781107378299,~0,~0.15),&\quad x<265,
    \end{cases}
  \end{equation*}
  and the state of the bubble is $$(\rho,u,v,p)=(0.1358,~0,~0,~0.05),\quad
  \sqrt{(x-45)^2+(y-45)^2}\leq 25.$$
  The setup of the second problem is the same except that the initial state of
  the bubble is replace with $(\rho,u,v,p)=(0.1358,~0,~0,~0.05)$.
\end{example}
Figure \ref{fig:SB1sch} shows the schlieren images of the
rest-mass density $\rho$ of the first shock-bubble interaction problem
at $t=90,180,270,360,450$ (from top to bottom) with $N_x=650,~N_y=180$.
Figure \ref{fig:SB2sch} gives the schlieren images of the
rest-mass density $\rho$ of the second shock-bubble interaction problem
at $t=100,200,300,400,500$ (from top to bottom) with $N_x=650,~N_y=180$.
Those plots clearly show the dynamics of the interaction between the
shock waves and the bubbles,
and the discontinuities and some small wave structures are also captured well by
our entropy stable scheme.

\begin{example}[Shock-vortex interaction]\label{ex:SV}\rm
  The final example is about the interaction between a shock wave and a vortex.
  The computational domain is $[-17,3]\times[-5,5]$, with reflective boundary
  conditions at $y=\pm5$, inflow and outflow boundary conditions at $x=3$ and
  $x=-17$, respectively, and the adiabatic index is $\Gamma=1.4$.
  We put a similar isentropic vortex initially centered at $(0,0)$ as in Example
  \ref{ex:acc2D}, except that the vortex here is moving left with magnitude of $w=0.9$.
  A planar stationary shock wave with Mach number $M_s=1.5$ is placed at $x=-6$,
  which is initially far away from the vortex, thus the pre-shock state is a
  const state $(\rho,u,v,p)=(1,-0.9,~0,~1)$.
  Then from the jump condition and the Lax shock condition, we can obtain the
  post-shock state as
  \begin{equation*}
    (\rho,u,v,p)=(4.891497310766981, -0.388882958251919,~0,~11.894863258311670).
  \end{equation*}
\end{example}

\begin{figure}[!ht]
  \centering
  \includegraphics[width=0.45\textwidth, trim=50 40 40 50, clip]{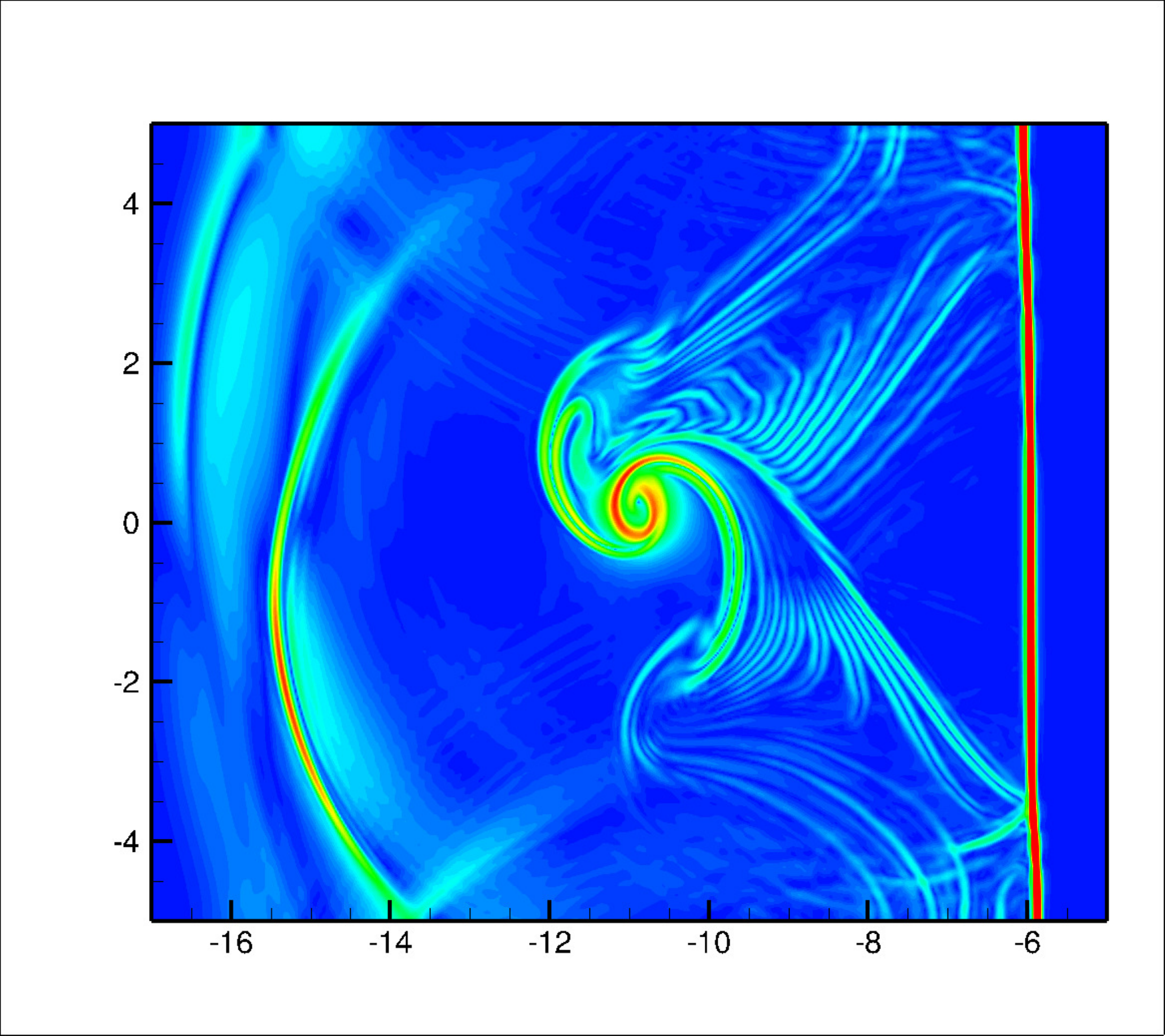}
  \includegraphics[width=0.45\textwidth, trim=50 40 40 50, clip]{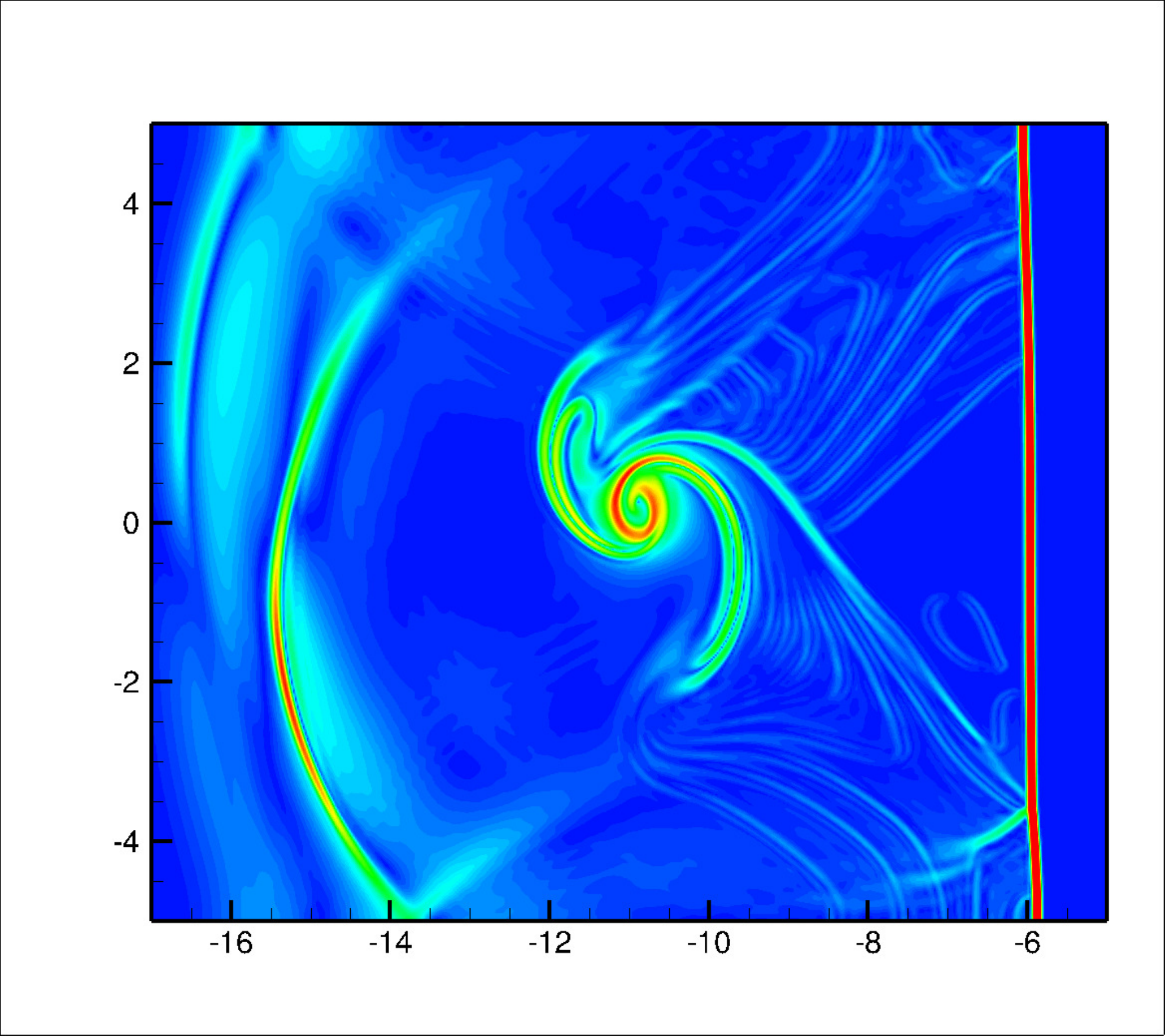}
  \caption{Example: \ref{ex:SV}: $50$ equally spaced contours lines from $0$ to
    $1$ of $\log_{10}(1+\abs{\nabla{\rho}})$ at $t=19$. $M_s=1.5$.
    $N_x=800,~N_y=400$.
  Left: entropy stable scheme; right: fifth order WENO-LLF.}
  \label{fig:SV}
\end{figure}

 Figure \ref{fig:SV} plots the contours with $50$ equally spaced contour
lines from $0$ to $1$ of $\log_{10}(1+\abs{\nabla{\rho}})$, obtained with our
entropy stable scheme and the fifth-order finite difference WENO scheme with
the local Lax-Friedrichs splitting. The computation is performed until $t=19$ with
$N_x=800,~N_y=400$. We can see that, after the interaction of the vortex and the
shock, the shock wave is still located at $x=-6$, and many linear and non-linear
waves, and sound waves generate and propagate in the domain.
Our entropy stable scheme can
capture the subtle details better than the fifth-order finite difference WENO scheme with
the local Lax-Friedrichs splitting.

\section{Conclusion}\label{section:conclusion}
%
For the special relativistic hydrodynamic (RHD) equations,
the schemes satisfying the discrete entropy
condition for a convex entropy function  has not been
considered before.
This paper has presented the  high-order accurate entropy stable finite difference schemes for
one- and two-dimensional special  RHD equations.
Those schemes
are built on the entropy conservative flux and the weighted essentially non-oscillatory (WENO)
technique as well as explicit Runge-Kutta time discretization.
The key is  to technically construct the {affordable} entropy conservative flux of the  semi-discrete second-order accurate entropy conservative schemes satisfying the semi-discrete entropy equality for the found convex entropy pair.
The entropy conservative schemes may become oscillatory near the shock wave, thus
as soon as  the entropy conservative flux  is derived,
the dissipation term can be added to give the semi-discrete entropy stable schemes satisfying the
semi-discrete entropy inequality with the given convex entropy function.
The WENO reconstruction  for the scaled entropy variables
and the high-order explicit Runge-Kutta time discretization
are implemented to obtain the fully-discrete high-order schemes.
%
%
Several numerical tests are conducted to validate the
accuracy and the ability to capture discontinuities  of our
entropy stable schemes. Especially, the Shock-vortex interaction problem
is designed for the first time.
The results show that
our schemes can achieve designed accuracy, and can well resolve the
discontinuities and subtle details.
In future, it will be interesting to study the physical-constraint-preserving
property of the entropy stable schemes, or extend them to the relativistic
magnetohydrodynamic equations.

\section*{Acknowledgments}
The authors were partially supported by the Special Project on High-performance Computing under the
National Key R\&D Program (No. 2016YFB0200603), Science Challenge Project (No. JCKY2016212A502),
the National Natural Science Foundation of China (Nos. 91630310 \& 11421101),
and High-performance Computing Platform of Peking University.


\end{document}